\documentclass[9pt,c,3p]{elsarticle}
\usepackage{amsmath}
\usepackage{amsfonts}
\usepackage{mathrsfs}
\usepackage{graphicx}
\usepackage{multirow}
\usepackage[english]{babel}
\usepackage{times}
\usepackage{enumerate}
\usepackage{a4wide}
\usepackage{amssymb}
\usepackage{amsbsy}
\usepackage{pgf}
\usepackage[textsize=footnotesize,color=yellow]{todonotes}
\usepackage{subfigure}
\usepackage{url}		
\usepackage{listings}
\usepackage{tikz}
\usepackage{pgfplots}
\usepackage[space]{grffile}
\usepackage{forloop}
\usepackage[hidelinks]{hyperref}

\newcommand{\vect}[1]{\ensuremath\boldsymbol{#1}}
\newcommand{\NVRtensor}[1]{\vect{#1}}

\newcommand{\norm}[1]{\left\Vert #1 \right\Vert}
\newcommand{\NVRgrad}{\nabla}
\newcommand{\NVRdiv}{\NVRgrad \cdot}

\newcommand{\eqdef}{\stackrel{\text{\tiny def}}{=}}

\newcommand{\NVRHdiv}{\ensuremath H(\text{\rm div})\,}

\definecolor{lightlightgray}{gray}{0.95}
\definecolor{lightlightblue}{rgb}{0.4,0.4,0.95}
\definecolor{lightlightgreen}{rgb}{0.8,1,0.8}
\newcommand{\eval}[2][\right]{\relax
  \ifx#1\right\relax \left.\fi#2#1\rvert}

\catcode`@=12

\let\tilde\widetilde

\begin{document}
\baselineskip=16pt
\parskip= 4pt

\newpage


\title{A Geometric Multigrid Preconditioning Strategy for DPG System Matrices}
\author[anl]{Nathan V. Roberts\corref{cor1}}
\author[rice]{Jesse Chan}
\cortext[cor1]{Corresponding author.\\  \emph{Email addresses}: nvroberts@anl.gov, Jesse.Chan@caam.rice.edu.}

\address[anl]{Argonne Leadership Computing Facility, Argonne National Laboratory, Argonne, IL, USA.}
\address[rice]{Rice University, Houston, TX, USA}

\begin{abstract}
The discontinuous Petrov-Galerkin (DPG) methodology of Demkowicz and Gopalakrishnan \cite{DPG1,DPG2} guarantees the optimality of the solution in an energy norm, and provides several features facilitating adaptive schemes.   A key question that has not yet been answered in general---though there are some results for Poisson, e.g.---is how best to precondition the DPG system matrix, so that iterative solvers may be used to allow solution of large-scale problems.

In this paper, we detail a strategy for preconditioning the DPG system matrix using geometric multigrid which we have implemented as part of Camellia \cite{Camellia}, and demonstrate through numerical experiments its effectiveness in the context of several variational formulations.  We observe that in some of our experiments, the behavior of the preconditioner is closely tied to the discrete test space enrichment.

We include experiments involving adaptive meshes with hanging nodes for lid-driven cavity flow, demonstrating that the preconditioners can be applied in the context of challenging problems.  We also include a scalability study demonstrating that the approach---and our implementation---scales well to many MPI ranks.
\end{abstract}

{\bf Key words:}
Discontinuous Petrov Galerkin, adaptive finite elements, iterative solvers, geometric multigrid

{\bf AMS subject classification:} 

\maketitle

\paragraph{Acknowledgments} This work was supported by the Office of Science, U.S. Department of Energy, under Contract DE-AC02-06CH11357.  This research used resources of the Argonne Leadership Computing Facility, which is a DOE Office of Science User Facility supported under Contract DE-AC02-06CH11357.

\section{Introduction}
\label{sec:intro}
The discontinuous Petrov-Galerkin (DPG) methodology of Demkowicz and Gopalakrishnan \cite{DPG1,DPG2} minimizes the residual of the solution in an energy norm, and has several features facilitating adaptive schemes.  For well-posed problems with sufficiently regular solutions, DPG can be shown to converge at optimal rates---the inf-sup constants governing the convergence are mesh-independent, and of the same order as those governing the continuous problem \cite{DPGStokes}. DPG also provides an accurate mechanism for evaluating the residual error which can be used to drive adaptive mesh refinements.

DPG has been studied for a host of PDE problems---including Poisson \cite{DPG6}, linear elasticity \cite{BramwellDemkowiczGopalakrishnanQiu11}, Stokes \cite{DPGStokes}, and compressible \cite{DPGCompressible} and incompressible Navier-Stokes \cite{DPGNavierStokes} problems, to name a few.  For each of these optimal convergence rates have either been proved a priori or observed in numerical experiments---in some cases, the solutions are even nearly optimal in terms of the absolute $L^2$ error (not merely the rate).  The global system matrices that arise from DPG formulations are symmetric (Hermitian) positive definite, making them good candidates for solution using the conjugate gradient (CG) method.  However, these matrices often have fairly large condition numbers which scale as $\frac{1}{h^{2}}$ (see \cite{GopalakrishnanQiu11} for the scaling estimate, and Table 9.3 in \cite{RobertsDissertation} for measurements), so that a good preconditioner is required before CG can be used effectively.

For us, a key motivation in the present work is the scalability of our solvers in Camellia \cite{Camellia}.  Prior to developing the preconditioners presented here, direct solvers were almost exclusively employed.  These solvers only scale to a certain limited system size, and can require substantially more memory than iterative solvers.  This is of particular concern for high-performance computing systems, where the memory per core is increasingly limited---for example, Argonne's Mira supercomputer has a BlueGene/Q architecture that has just one gigabyte per core.

The structure of the paper is as follows.  In Section \ref{sec:litReview}, we review some previous work in preconditioning DPG and similar systems.  In Section \ref{sec:dpgIntro}, we briefly introduce the DPG methodology and state the variational formulations we use for our numerical experiments.  In Section \ref{sec:gmgDetails}, we detail our geometric multigrid preconditioners.  In Section \ref{sec:numericalExperiments}, we present a wide variety of numerical experiments demonstrating the effectiveness of the approach.  We examine the scalability of our implementation in Section \ref{sec:scalability}.  We conclude in Section \ref{sec:conclusion}.  Some notes on our implementation in Camellia can be found in \ref{sec:implementation}; for reference, we provide numerical values for the smoother weights we employ in \ref{sec:schwarzWeights}.

\section{Literature Review}
\label{sec:litReview}

DPG methods incorporate aspects of both least squares and substructured finite element methods.  For least squares finite element methods \cite{bochev2009least}, multigrid methods have been the preconditioner of choice.  These were popularized as \textit{black-box} solvers for First Order Least Squares (FOSLS) finite element methods, as the elliptic nature of least-squares finite element formulations implies optimal convergence for both additive and multiplicative multigrid methods for second order partial differential equations \cite{cai1997first}.  For finite element discretizations based on the mesh skeleton, such as static condensation, mortar, or hybridized methods, Schur complement or substructuring preconditioners based on the mesh skeleton have been developed \cite{achdou1999iterative}.  For a comprehensive review of such preconditioners, we direct the interested reader to \cite{toselli2005domain,BrennerScott}.  

Wieners and Wohlmuth examine preconditioning of the substructured DPG system \cite{WienersWohlmuth} --- this is equivalent to the Schur complement/static condensation we will discuss in Section \ref{sec:gmgDetails} --- and prove that given an effective preconditioner of the original DPG system, an effective preconditioner for the Schur complement system can be derived.  This construction involves three ingredients: a trace operator which extracts boundary traces of functions on the mesh skeleton, a secondary space whose image under the trace operator yields the DPG trace space, and a self-adjoint preconditioner for the secondary space.  Since the dual of the trace operator maps traces to the secondary space, a preconditioner for the Schur complement system can be applied to trace unknowns by extending them to the secondary space (using the dual of the trace operator), applying the preconditioner on the secondary space, and applying the trace operator to map the results back to the trace space.  

The present work differs from the aforementioned literature in that a multi-level preconditioning is directly applied to the Schur complement system for the trace unknowns.  The smoother is an overlapping additive Schwarz domain decomposition method, for which mesh and order independence has been shown for the Poisson equation under a fixed subdomain overlap width \cite{BarkerDPG} .  
Fischer and Lottes use a multigrid preconditioner with Schwarz smoothing for a fractional-step Navier-Stokes solver, in which the two steps involve Poisson solves \cite{Fischer2005}.  Our approach follows theirs in several respects, though notably we omit a weighting matrix that they introduce (the one they refer to as $W$), as we found it to be detrimental in the context of our solvers.  

We note that the preconditioning strategy presented here is \textit{black-box}, in the sense that it can be applied to any DPG system matrix.  However, we observe that the performance of this preconditioner worsens for singularly perturbed differential equations, which often require more specialized techniques.  Examples of singularly perturbed differential equations include the frequency-domain Helmholtz equation or convection-diffusion equation with small viscosity.   The preconditioning of Helmholtz equations is addressed by Gopalakrishnan and Schoberl \cite{gopalakrishnan2015degree}, who observe wavenumber and $p$-independence on a fixed grid for a one-level multiplicative Schwarz preconditioner involving forward-backward Gauss-Seidel sweeps.  Similar wavenumber independence is observed by Li and Xu \cite{xu2016domain} for a one-level additive preconditioner.  The preconditioning of DPG for convection-diffusion problems is an open problem, and will likely involve streamline-aware techniques \cite{bey1997downwind,loghin1997preconditioning}, though these may be simplified by the self-adjoint nature of the DPG least squares discretization.  

\section{DPG with Ultraweak Variational Formulations}
\label{sec:dpgIntro}
In this section, we first provide a brief review of the DPG formulation for the Poisson problem.  We then turn to defining the DPG method for an abstract variational formulation, and defining the \emph{graph norm} on the test space in the abstract setting.  Finally, we state the formulations for Poisson, Stokes, and Navier-Stokes that we employ in the numerical experiments that follow.  Here, we aim simply to specify the operational approach; for full details of the functional settings, we refer the reader to previous treatments of the Poisson \cite{DPG6,GopalakrishnanQiu11} and Stokes formulations that we employ here \cite{DPGStokes}.

\subsection{The Ultraweak Variational Formulation for the Poisson Problem}
\label{sec:PoissonFormulation}
Consider the problem
\begin{align*}
-\Delta u &= f \quad {\rm in } \quad \Omega, \\
u &= u_g  \quad {\rm on } \quad \partial\Omega.
\end{align*}
First, we rewrite as a first-order system by introducing $\sigma = \nabla u$, giving us:
\begin{align*}
-\NVRdiv \sigma &= f \quad {\rm in } \quad \Omega, \\
\sigma - \nabla u &= 0 \quad {\rm in } \quad \Omega, \\
u &= u_g  \quad {\rm on } \quad \partial\Omega.
\end{align*}
Suppose some finite element mesh $\Omega_h$ is given.  We then multiply the strong equations by test functions $v$ and $\tau$ and integrate by parts element-wise to get:
\begin{align*}
(\sigma, \nabla v)_{\Omega_h} - \langle \sigma \cdot \vect{n}, v \rangle_{\partial \Omega_h}  &= (f, v)_{\Omega_h}, \\
(\sigma, \tau)_{\Omega_h} + (u, \NVRdiv \tau)_{\Omega_h} - \langle u, \tau \cdot \vect{n} \rangle_{\partial \Omega_h} &= 0.
\end{align*}
To satisfy regularity requirements such that we may place $u \in L^2(\Omega), \sigma \in \vect{L}^2(\Omega)$, we introduce new unknowns $\widehat{u}, \widehat{\sigma}_n$ on the mesh skeleton $\partial \Omega_h$.  Summing the equations, we have
\begin{align*}
b((u,\sigma,\widehat{u}, \widehat{\sigma}_n), (v, \tau))
&\eqdef (\sigma, \nabla v)_{\Omega_h} - \langle \widehat{\sigma}_n, v \rangle_{\partial \Omega_h} +  (\sigma, \tau)_{\Omega_h} + (u, \NVRdiv \tau)_{\Omega_h} - \langle \widehat{u}, \tau \cdot \vect{n} \rangle_{\partial \Omega_h}\\
&= (f, v)_{\Omega_h}.
\end{align*}
$b(\cdot,\cdot)$ is then referred to as the \emph{ultraweak variational formulation}: all differential operators have been moved to the test space through integration by parts, and new unknowns (known as \emph{trace variables}) have been introduced on the mesh skeleton, resulting in an energy setting wherein the variables defined on the volume---the \emph{field variables}---lie in $L^2$ spaces.


\subsection{The DPG Method for an Abstract Variational Problem}
Suppose now that we have some variational problem $b(u,v) = l(v)$ where $u \in U, v \in V$ for $U, V$ (infinite-dimensional) Hilbert spaces.  Suppose further that some discretization $U_h \subset U$ of the trial space is given on a finite element mesh $\Omega_h$.  The space $V$ is ``broken'' in the sense that test functions $v$ are only required to be conforming element-wise; they are allowed to be discontinuous across elements in $\Omega_h$.  Let $(\cdot, \cdot)_V$ be the inner product on the test space.  The \emph{ideal DPG method} consists of solving, for every basis function $e \in U_h$, the problem
\begin{align*}
(v_{e}, v)_V = b(e, v) \quad \forall v \in V,
\end{align*}
and using the solutions $v_e$ as the discrete test space for the problem.  The test functions $v_e$ are exactly the Riesz representations of the bilinear form $b(e, \cdot)$, interpreted as a functional on the test space.  At this point, $V$ is still infinite-dimensional, and therefore solving for $v_e$ is impractical.  If we introduce a discrete space $\tilde{V}_h \subset V$, then the problem
\begin{align*}
(\tilde{v}_{e}, v)_V = b(e, v) \quad \forall v \in \tilde{V}_h
\end{align*}
may be solved discretely and element-wise.  Using $V_h \eqdef \left\{\tilde{v}_e : e \in U_h \right\}$ as the discrete test space for the problem, we arrive at the \emph{practical DPG method} \cite{GopalakrishnanQiu11}.  For this to work well, the ``practical'' test functions $\tilde{v}_e$ should approximate the ``ideal'' $v_e$.  To achieve this in practice, we enrich the polynomial order of the trial space by some $\Delta k$, and use this as the polynomial order for the discrete test space $\tilde{V}_h$---see Section \ref{sec:polynomialOrder} for a precise definition of $\Delta k$.  The appropriate choice for $\Delta k$ is problem-dependent; following previous work on Poisson and linear elasticity \cite{GopalakrishnanQiu11}, a reasonable starting point is the spatial dimension $d$.  However, as we will see in our Navier-Stokes experiments, this may have consequences for the behavior of the preconditioners that are the subject of the present work.

\subsection{The Graph Norm}
In the above, we have left unspecified what the inner product on the test space should be.  DPG minimizes the residual in the energy norm
\begin{align*}
\norm{u}_E \eqdef \norm{b(u,\cdot)}_{V^{'}} = \sup_{v \ne 0} \frac{b(u,v)}{\norm{v}_V};
\end{align*}
the norm on $V$ thus determines the norm in which the residual is minimized.  If we want to minimize the error in the $L^2$ norm of the field variables, the \emph{graph norm} on the test space $V$ is a good choice for many problems.  Suppose the original first-order system is of the form $L u = f$.  After we multiply by test variable $v$, integrate by parts, and introduce trace variables, the weak system takes the form
\begin{align*}
b(u,v) = (u, L^{*} v) + \{\text{boundary terms}\},
\end{align*}
where $L^{*}$ is the adjoint of $L$.  The graph norm of $v$ is then defined by
\begin{align*}
\norm{v}_{\rm graph}^2 = \norm{L^{*} v}^2 + \beta \norm{v}^2,
\end{align*}
where $\beta$ is a scaling parameter.  In all our numerical experiments, we use the graph norm on the test space with $\beta = 1$.  For full details, including a proof that using the graph norm on the test space suffices to guarantee optimal $L^2$ convergence rates for a wide class of problems, see Roberts et al. \cite{DPGStokes}.

\subsection{The Polynomial Order $k$}
\label{sec:polynomialOrder}
Throughout this paper, the polynomial order $k$ of a mesh refers to the order of the field variables.  In the Poisson formulation above, we define two trace variables, $\widehat{u}$---the trace of an $H^1$ variable---and $\widehat{\sigma}_n$---the normal trace of an $\NVRHdiv$ variable.  When field variables take polynomial order $k$, the $\NVRHdiv$ traces will then also have order $k$, while the $H^1$ trace variables will have order $k+1$.  We select the polynomial orders in this way so that all variables will converge in $L^2$ at the same rate---for full details, see \cite[Section 3.1]{DPGStokes}.  A lowest-order mesh will have constant field variables.  The test space enrichment $\Delta k$ is taken relative to the $H^1$ order, so that an $H^1$ test variable will have order $k+1+\Delta k$.

\subsection{Variational Formulations}
Here, we specify the variational formulations we use in our experiments.
\paragraph{Poisson}
Our variational formulation for Poisson is as specified above:
\begin{align*}
b((u,\sigma,\widehat{u}, \widehat{\sigma}_n), (v, \tau))
&\eqdef (\sigma, \nabla v)_{\Omega_h} - \langle \widehat{\sigma}_n, v \rangle_{\partial \Omega_h} +  (\sigma, \tau)_{\Omega_h} + (u, \NVRdiv \tau)_{\Omega_h} - \langle \widehat{u}, \tau \cdot \vect{n} \rangle_{\partial \Omega_h}\\
&= (f, v)_{\Omega_h}.
\end{align*}

\paragraph{Stokes}
For Stokes, we use the velocity-gradient-pressure formulation specified in \cite{DPGStokes}:
\begin{align*}
b_{\rm Stokes}(u,v) \eqdef &\left(\NVRtensor{\sigma} - p \NVRtensor{I}, \NVRgrad \vect{v} \right)_{\Omega_{h}} - \left\langle \widehat{\vect{t}}_{n}, \vect{v} \right\rangle_{\partial \Omega_{h}}\\
&+\left(\vect{u}, \NVRgrad q \right)_{\Omega_{h}} -  \left\langle \widehat{\vect{u}} \cdot \vect{n}, q  \right\rangle_{\partial \Omega_{h}}\\
&+\left( \NVRtensor{\sigma}, \NVRtensor{\tau} \right)_{\Omega_{h}} + \left( \mu \vect{u}, \NVRdiv \NVRtensor{\tau} \right)_{\Omega_{h}} - \left\langle \widehat{\vect{u}}, \NVRtensor{\tau} \vect{n}  \right\rangle_{\partial \Omega_{h}} =  (\vect{f}, \vect{v} )_{\Omega_{h}},
\end{align*}
where we define group variables $u = (\vect{u},\NVRtensor{\sigma},p,\widehat{\vect{u}},\widehat{\vect{t}}_n)$ and $v = (\vect{v},\NVRtensor{\tau},q)$, and $\mu$ is the (constant) viscosity, $\vect{u}$ is the velocity, $p$ is the pressure, $\NVRtensor{\sigma}$ is the gradient of the $\mu$-weighted velocity, $\widehat{\vect{u}}$ is the velocity trace, and $\widehat{\vect{t}}_n$ is a pseudo-traction, the trace of $( \NVRtensor{\sigma} - p \NVRtensor{I} ) \vect{n}$.

\paragraph{Navier-Stokes}
For Navier-Stokes, we also use a velocity-gradient pressure formulation, based on the Stokes formulation (for Reynolds number ${\rm Re}$, we take viscosity $\mu = \frac{1}{\rm Re}$) and used in \cite{DPGNavierStokes}.  The nonlinear formulation is given by
\begin{align*}
b_{\rm Stokes}(u,v) + {\rm Re \,} (\vect{u} \cdot \NVRtensor{\sigma}, \vect{v}) = (\vect{f}, \vect{v}).
\end{align*}
Linearizing about $u + \Delta u$, we have:
\begin{align*}
b_{\rm Stokes}(\Delta u,v) + {\rm Re \,} (\Delta \vect{u} \cdot \NVRtensor{\sigma} + \vect{u} \cdot \Delta \NVRtensor{\sigma}, \vect{v})_{\Omega_{h}} &= (\vect{f}, \vect{v})_{\Omega_{h}} - b_{\rm Stokes}(u,v) - {\rm Re \,} (\vect{u} \cdot \NVRtensor{\sigma}, \vect{v})_{\Omega_{h}}.
\end{align*}
We solve the original nonlinear problem by a standard Newton iteration; given an initial guess $u = u_0$, we iterate by solving the linearized problem for increment $\Delta u$ and setting $u := u + \Delta u$, continuing until some stopping criterion on $\Delta u$ is met.

\section{Our Multigrid Operators}
\label{sec:gmgDetails}
Because DPG always results in a symmetric (Hermitian) positive definite system matrix, we may use the conjugate gradient method to solve the global system iteratively.  To do so efficiently, however, requires a good preconditioner.  In Camellia, we aim to provide implementations that are applicable across a wide range of PDE problems---essentially, we want to provide good defaults that the user may override when he or she has special requirements.

Because of their broad applicability (especially for the class of least squares finite element methods to which DPG belongs), and because our adapted meshes include a natural hierarchy of geometric refinements, we use geometric multigrid preconditioners for conjugate gradient solves.  Our numerical experiments suggest that this approach works well in the context of many problems; we have used this approach with Poisson, Stokes, linear elasticity, and both compressible and incompressible Navier-Stokes.

Specification of a multigrid operator involves the following choices:
\begin{itemize}
\item a prolongation operator $P$ to transfer data from the coarse mesh to the fine,
\item a restriction operator $R$ to transfer data from the fine mesh to the coarse,
\item a smoother $B$ to operate locally on the fine mesh, and
\item a multigrid strategy that specifies the way that the various operators work together.
\end{itemize}

For the remainder of this work, we take $R=P^T$.  If the fine-grid system matrix is $A$, then the coarse-grid system matrix is $A_C = P^T A P$.

The multigrid strategy specifies both the way that the hierarchy is traversed in the course of one multigrid iteration and the way that at each level the result of the coarse solve is combined with the result of the smoother.  Hierarchy traversal choices include V-cycle, W-cycle, and full multigrid; the standard combination choices are multiplicative and additive.  Additive combinations (also known as two-level) apply both smoother and coarse operator to the same residual, while multiplicative combinations apply them in sequence, recomputing the residual in between.  Additive combinations have the advantage of avoiding recomputation of the residual, typically with the tradeoff that more iterations are required to achieve a specified tolerance.

Our preferred multigrid strategy is the multiplicative V-cycle,\footnote{We have found that full multigrid and W-cycles do reduce iteration counts, but with costs per iteration which are substantially higher.  Using an additive V-cycle, on the other hand, gives us a preconditioner that does not scale---as the grid size $h$ decreases or the polynomial order $k$ increases, the iteration counts required to achieve a specified residual tolerance grow without bound.} which for iteration $n$ with current solution $x_{n-1}$ proceeds as follows:
\begin{enumerate}
\item Compute the residual: $r_n^0 = b - A x_{n-1}$.
\item Apply the smoother: $x_n^1 = x_{n-1} + \sigma B r_n^0$.
\item Recompute the residual: $r_n^1 = b - A x_n^1$.
\item Apply the coarse operator: $x_n^2 = x_n^1 + P M_C P^T r_n^1$.
\item Recompute the residual: $r_n^2 = b - A x_n^2$.
\item Apply the smoother: $x_n = x_n^3 = x_n^2 + \sigma B r_n^2$.
\end{enumerate}

At the coarsest level, $M_C = A_C^{-1}$ is a direct solve; at intermediate levels, $M_C$ will be the application of steps 1-6 with the operators corresponding to the next-coarsest grid.  It remains to specify the prolongation operator $P$, the smoother $B$, and the smoother weight $\sigma$.  We describe our choices for each of these in turn.

\subsection{The Prolongation Operator $P$}
A natural, minimal requirement for the prolongation operator $P$ is that a solution which is exact on the coarse mesh should, when prolongated to the fine mesh, remain an exact solution.  In the usual finite element case with nested discrete spaces, this is straightforward for both $h$ and $p$ refinements: each variable in the fine mesh is also defined on the coarse, and every coarse basis function may be represented exactly in terms of fine basis functions.  The prolongation operator is then the one that takes a coarse basis function to its representation in terms of fine basis functions.

In the case of DPG for $p$ refinements, the same facts hold.  However, for $h$ refinements, we have an additional complication: the introduction of new faces implies the introduction of new trace variables on the fine mesh, which have no predecessor on the coarse mesh.  Simply using the prolongation operator as defined in the usual case would violate our exact solution requirement: the new trace variables on the fine mesh would be initialized to zero, resulting in a nonzero residual.  Now, each trace variable, as the term implies, is a trace of some combination of field variables.\footnote{It is typically the case that the variables traced have higher regularity requirements than those of the field variables employed in the actual DPG computation.  In the case of ultraweak formulations, for example, each field variable requires only $L^2$ regularity, of which we may not take a trace.  We do not yet have theoretical justification for this formal violation of regularity requirements, though we would note that in the discrete setting all the basis functions are polynomials, so that the trace operator is well-defined.}

For each interior face in an $h$-refinement, we may therefore define an operator $\gamma(U)$ that maps the field variables $U$ to the traces $\widehat{U}$.  Wherever new variables are defined on the interior of a coarse-grid element, we use $\gamma$ to prolongate from the field bases on the interior of the coarse-grid element to the trace variables on the newly-defined interfaces of the fine grid.\footnote{It is worth emphasizing that we only use $\gamma$ to prolongate field degrees of freedom to the trace degrees of freedom that lie strictly in the interior of the coarse element---when an $H^1$ trace lies on the interface between an interior face and the exterior of the coarse element, for example, the coarse element's trace degrees of freedom suffice for prolongation.}  This allows us to satisfy our exact solution requirement.

\subsection{Static Condensation}
We have found it beneficial to employ \emph{static condensation}, which reduces the size of the global system by locally eliminating the interior degrees of freedom on each element (for ultra-weak DPG, these are precisely the field degrees of freedom).  Suppose that our discrete DPG system is of the form $Kx = F$.  The matrix can be reordered to take the form
\[
\left(\begin{array}{cc}
K_{11} &K_{12}\\
K_{12}^T &K_{22}
\end{array}
\right)
\left(\begin{array}{c}
u\\
f
\end{array}
\right) = 
\left(\begin{array}{c}
F_{1}\\
F_{2}
\end{array}
\right)
\]
where $K_{11}$ is block diagonal, $u$ represents the degrees of freedom corresponding to field variables, and $f$ represents those corresponding to the trace variables. Noting that $u = K_{11}^{-1}(F_{1}-K_{12}f)$, we can substitute this into the equation $K_{12}^{T}u + K_{22}f = F_{2}$ to obtain the Schur complement system for the trace degrees of freedom:
\[
(K_{22}-K_{12}^TK_{11}^{-1}K_{12})f = F_{2}-K_{12}^TK_{11}^{-1}F_{1}.
\]
Since $K_{11}$ is block-diagonal, its inversion can be carried out element-wise and in parallel; usually, $K_{12}$ is a significantly smaller matrix and the computational cost of the global solve is reduced.

The principal benefit of static condensation is that the operators $P$ and $B$ become less expensive to store and to apply.  The computational cost to determine $P$ in the context of $h$-refinements does increase, however, because before we may apply the coarse field to fine trace operator $\gamma$ we must first compute and apply the coarse trace to coarse field operator, which is given by the static condensation formula $u = - K_{11}^{-1} K_{12} f$.  (Note that in an iterative context, $F_1$ may be neglected.)

\subsection{The Smoother $B$}
For both $h$ and $p$-multigrid, we employ an overlapping additive Schwarz preconditioner.  In the context of our two-grid experiments, we have found that for many problems a minimal-overlap operator scales for $p$-multigrid, while for $h$-multigrid a 1-overlap Schwarz operator is required.  By minimal overlap, we mean that each Schwarz block corresponds to the degrees of freedom seen by an element, including those shared with its neighbors (the overlap region is therefore the element boundary).  By 1-overlap, we mean that each Schwarz block corresponds to the degrees of freedom seen by an element and its face neighbors.\footnote{When a face contains a hanging node, the face neighbors of the coarse element are the fine elements which have faces resulting from the refinement of the coarse face.}

\subsection{The $\sigma$ Weight: Bounding Eigenvalues of the Schwarz-Preconditioned System}\label{sec:sigmaWeightDiscussion}
We have found it important to weight our Schwarz-preconditioned matrix, $BA$, with a weight $\sigma$ such that the maximum eigenvalue of $\sigma B A$ is at most 1.  (If we do not do this, our preconditioner is no longer guaranteed to be positive definite.)  Below, we describe first a conservative choice for $\sigma$---one for which we can prove the bound in general---and a more aggressive choice that we have found to work well in practice.

Smith et al. \cite{SmithBjorstadGropp} have shown\footnote{The main result, a bound for the eigenvalue in a general context, is Smith et al.'s Lemma 3 (p. 156); the particular value we use here is derived from their coloring argument found on p. 167.} that the largest eigenvalue of $B A$ is bounded above by $N_c + 2$, where $N_c$ is the number of colors required to color the Schwarz domains such that no two adjacent domains have like color.  If we consider the adjacency graph $G$ for the Schwarz domains, it is clear that if we count the number $N_i$ of Schwarz blocks in which degree of freedom $i$ participates, then $N_{\rm max} \eqdef \max_i N_i = \Delta(G)$, where $\Delta(G)$ is the maximum degree in the graph.  Supposing that $G$ is neither complete nor an odd cycle---which is almost always true in practice---then $N_c \leq \Delta(G)$ (this is Brooks's Theorem \cite{Brooks1941}).  We then may take $\sigma = \frac{1}{N_{\rm max} + 2}$ to bound the eigenvalue of $\sigma B A$ as desired.

It is worth noting that the bound on the eigenvalue may be loose, and the above $\sigma$ may therefore be suboptimal.  In the present work, we employ a more aggressive choice for $\sigma$ that is computed as follows.  For each Schwarz domain, count the \emph{face neighbors} of cells in that domain (including the cells that belong to the domain).  Call the maximum such count $N$; we then take $\sigma$ to be $\frac{1}{N+1}$.  For reference, values of $\sigma$ on uniform grids are provided in \ref{sec:schwarzWeights}.

\subsection{Determining a Mesh Hierarchy}
\label{sec:meshHierarchy}
Given a particular fine mesh with polynomial order $k_{\rm fine}$, how do we determine an appropriate mesh hierarchy?

As noted above, in our two-grid experiments, we have found that for many problems minimal-overlap additive Schwarz scales for $p$-multigrid but not for $h$-multigrid (for $h$-multigrid typically 1-overlap additive Schwarz is required).  Because 1-overlap additive Schwarz can be fairly expensive for higher-order meshes, we therefore prefer a mesh hierarchy that performs $p$-coarsening first, followed by $h$-coarsenings.  This allows us to limit our application of 1-overlap operators to lower-order meshes.

Here, we describe Camellia's default approach for generating a mesh hierarchy from a provided fine mesh---we produce a stack of meshes, from coarse to fine.  We usually take the polynomial order on the coarsest mesh, $k_{\rm coarse}$, to be $1$ or $0$.  We also define a boolean parameter, \texttt{skipIntermediateP}, which governs whether more than one $p$-coarsening in performed.  (We recommend choosing \texttt{skipIntermediateP = true}; while iteration counts are generally a bit higher, the computational cost of determining, storing, and applying the operator is reduced.)

\begin{enumerate}
\item Let the coarsest mesh be the set of original, unrefined cells (the ``root'' mesh geometry), with polynomial order $k = k_{\rm coarse}$.  Add it to the stack.
\item While there are elements in the last mesh on the stack coarser in $h$ than the fine mesh:
    \begin{enumerate}
    \item Duplicate the last mesh.
    \item $h$-refine (once) any elements that are coarser in $h$ than the fine mesh.
    \item Add the resulting mesh to the stack.
    \end{enumerate}
\item If \texttt{skipIntermediateP} is \texttt{true}: add the fine mesh to the stack.
\item If \texttt{skipIntermediateP} is \texttt{false}: while there are elements in the last mesh on the stack coarser in $p$ than the fine mesh:
    \begin{enumerate}
    \item Duplicate the last mesh.
    \item For any element for which $k < k_{\rm fine}$, $p$-refine to $\min(2 * k, k_{\rm fine})$.
    \item Add the resulting mesh to the stack.
    \end{enumerate}
\end{enumerate}
An example mesh hierarchy can be found in the context of our cavity flow experiments, in Figure \ref{fig:cavityFlowHierarchy}.

\subsection{The CG Stopping Criterion}
As is standard practice, in all our experiments we employ a stopping criterion based on the $\ell_2$-norm of the discrete residual vector $r$, scaled by the $\ell_2$-norm of the discrete right-hand side, $b$.  However, it is worth noting that this choice may not be optimal in terms of minimizing the error in the discrete energy norm; Arioli has proposed an alternative stopping criterion \cite{Arioli2003}, which we are considering adopting in future work.

\section{Numerical Experiments}\label{sec:NumericalExperiments}
\label{sec:numericalExperiments}
We present a series of numerical experiments with the multigrid-preconditioned conjugate gradient solver described above and implemented in Camellia.  For each of Poisson, Stokes, and Navier-Stokes, we consider smooth non-polynomial solutions.  To gain some insight into the behavior of the individual operators, we begin by using exactly two grids.  We then turn to multigrid experiments, first with uniform refinements, then with adaptive refinements.  All experiments are carried out using hypercube meshes---quadrilateral meshes in two dimensions, and hexahedral meshes in three dimensions.  For both trial and test space basis functions, we use the conforming nodal bases provided by the Intrepid package \cite{Intrepid}, with nodes defined at the Lobatto points.

\subsection{Problem Definitions}
\label{sec:problemDefinitions}
In almost all the experiments that follow, we use problems with non-polynomial smooth solutions.  The one exception is in the adaptive refinements, where we use a lid-driven cavity flow problem.  The problems we use are described below.

\paragraph{Poisson}
For the Poisson problem, we use homogeneous boundary conditions $\widehat{\phi}|_{\partial \Omega} = 0$ and unit forcing $f = 1$ on the domain $\Omega=[0,1]^d$, where $d=1,2,$ or $3$ is the spatial dimension.

\paragraph{Stokes}
For the analytic Stokes solution, in three dimensions we employ the manufactured solution
\begin{align*}
u_{1} &= - e^x \, y \, \cos y + \sin y\\
u_{2} &= e^x \, y \, \sin y + e^z \, y \, \cos y\\
u_{3} &= - e^z (\cos y - y \sin y)\\
p &= 2 \mu \, e^{x} \sin y
\end{align*}
on domain $\Omega = (-1,1)^3$, taking $\mu=1$.  This solution is an extension of the two-dimensional manufactured solution employed by Cockburn et al. \cite{CockburnKanschatSchotzauSchwab03}, which we have also previously used \cite{DPGStokes}; the two-dimensional version can be arrived at by dropping the terms involving the $z$ coordinate.  For our two-dimensional experiments, this is the solution we use.  For the pressure, in lieu of a zero-mean constraint, we impose the discrete condition that $p=0$ at the origin (which is also the center of the domain).\\

\paragraph{Navier-Stokes}
For the analytic Navier-Stokes solution, we use the classical two-dimensional solution due to Kovasznay \cite{Kovasznay}:
\begin{align*}
u_{1} &= 1 - e^{ \lambda x } \cos (2 \pi y)\\
u_{2} &= \frac{\lambda}{2 \pi} e^{ \lambda x } \sin ( 2 \pi y )\\
p &= \frac{1}{2} e^{ 2 \lambda x } + C
\end{align*}
where $\lambda = \frac{\rm Re}{2} - \sqrt{ \left(\frac{\rm Re}{2}\right)^2 + (2 \pi)^2 }$.  We use $\Omega = (-0.5,1.5) \times (0,2)$ as our domain, and choose the constant $C$ to agree with the discrete constraint on the pressure (here, that it is zero at (0.5,1)).  As is common when studying Kovasznay flow, we use ${\rm Re}=40$.\\

\paragraph{Lid-Driven Cavity Flow for Stokes and Navier-Stokes}
For the experiments involving adaptivity, we use the Stokes and Navier-Stokes lid-driven cavity flow problem in two dimensions.  Details of the problem setup can be found in Section \ref{sec:adaptiveCavityFlow}.

\subsection{Two-Grid Experiments}
To investigate the behavior at each level of multigrid, we begin by examining the simplest case of two-level multigrid for $h$ and $p$.  For the $h$ case, we coarsen the fine grid once to produce a coarse grid.  For the $p$ case, we use a field polynomial order $k_{\rm coarse} = \lfloor k_{\rm fine} / 2 \rfloor$.  We use the smooth non-polynomial solutions described above for Poisson, Stokes, and Navier-Stokes.  We perform conjugate gradient iterations until we reach a tolerance of $10^{-10}$, as measured in the discrete $\ell_2$ norm, relative to the discrete problem's right hand side.

The results for Poisson in one, two, and three space dimensions can be found in Tables \ref{table:PoissonPTwoGridResults} and \ref{table:PoissonHTwoGridResults}.  In every case, the $p$-multigrid and $h$-multigrid operators \emph{scale}: as $k$ increases or $h$ decreases, the number of iterations required do not go up.

\begin{table}[h!b!p!]
\begin{center}
\begin{tabular}{| c  c  c  | c  c  c | c  c  c |}
\cline{1-9}
\multicolumn{9}{| c |}{$d=1$} \\
\cline{1-9}
\multicolumn{3}{| c |}{}	&$k$ & Mesh Width & Iterations &\multicolumn{3}{c |}{}\\
\multicolumn{3}{| c |}{}	&0,1,2,4,8,16	& 2,4,8,16,32,64	&1	&\multicolumn{3}{c |}{}\\
\cline{1-9}
\multicolumn{9}{| c |}{$d=2$} \\
\hline 
$k$ & Mesh Width & Iterations	&$k$ & Mesh Width & Iterations	&$k$ & Mesh Width & Iterations\\
\hline
1				&2	&4	&2				&2	&4	&4				&2	&6\\
				&4	&11	&				&4	&10	&				&4	&13\\
				&8	&17	&				&8	&13	&				&8	&14\\
				&16	&18	&				&16	&13	&				&16	&13\\
				&32	&18	&				&32	&12	&				&32	&13\\
				&64	&16	&				&64	&12	&				&64	&12\\
\hline
\multicolumn{9}{| c |}{$d=3$} \\
\cline{1-9}
$k$ & Mesh Width & Iterations	&$k$ & Mesh Width & Iterations	&\multicolumn{3}{c |}{}\\
\hline
1				&2	&8	&2				&2	&9	&\multicolumn{3}{c |}{}\\
				&4	&21	&				&4	&17	&\multicolumn{3}{c |}{}\\
				&8	&24	&				&8	&18	&\multicolumn{3}{c |}{}\\
				&16	&24	&				&16	&17	&\multicolumn{3}{c |}{}\\
				&32	&23	&				&32	&16	&\multicolumn{3}{c |}{}\\
\cline{1-9}
\end{tabular}
\end{center} 
\caption[]{Poisson $p$-multigrid with static condensation, two-grid results: iteration counts to achieve a residual tolerance of $10^{-10}$.  In each case, the coarse mesh is geometrically identical to the fine, but with polynomial order half that of the fine mesh ($0$ when the fine mesh has polynomial order $1$).  In every choice of polynomial order and mesh width, the 1D preconditioner effectively operates as an exact inverse, achieving the tolerance in a single CG iteration.}
\label{table:PoissonPTwoGridResults}
\end{table}

\begin{table}[h!b!p!]
\begin{center}
\begin{tabular}{| c  c  c  | c  c  c | c  c  c |}
\cline{1-9}
\multicolumn{9}{| c |}{$d=1$} \\
\cline{1-9}
$k$ & Mesh Width & Iterations	&$k$ & Mesh Width & Iterations &\multicolumn{3}{c |}{}\\
1	&2	&2	&2,4,8,16	&2	&1 &\multicolumn{3}{c |}{}\\
	&4	&3	&	&4	&3 &\multicolumn{3}{c |}{}\\
	&8	&5	&	&8	&5 &\multicolumn{3}{c |}{}\\
	&16	&6	&	&16	&5 &\multicolumn{3}{c |}{}\\
	&32	&7	&	&32	&7 &\multicolumn{3}{c |}{}\\
	&64	&7	&	&64	&6 &\multicolumn{3}{c |}{}\\
\cline{1-9}
\multicolumn{9}{| c |}{$d=2$} \\
\hline 
$k$ & Mesh Width & Iterations	& $k$ & Mesh Width & Iterations	& $k$ & Mesh Width & Iterations	\\
1	&2	&5	&2	&2	&5	&4	&2	&5\\
	&4	&12	&	&4	&13	&	&4	&14\\
	&8	&16	&	&8	&15	&	&8	&15\\
	&16	&16	&	&16	&14	&	&16	&15\\
	&32	&16	&	&32	&14	&	&32	&14\\
	&64	&16	&	&64	&13	&	&64	&14\\
\cline{1-9}
\multicolumn{9}{| c |}{$d=3$} \\
\hline 
$k$ & Mesh Width & Iterations	& $k$ & Mesh Width & Iterations	&\multicolumn{3}{c |}{}\\
1	&2	&7	&2	&2	&8	&\multicolumn{3}{c |}{}\\
	&4	&18	&	&4	&18	&\multicolumn{3}{c |}{}\\
	&8	&20	&	&8	&19	&\multicolumn{3}{c |}{}\\
	&16	&21	&	&16	&19	&\multicolumn{3}{c |}{}\\
	&32	&21	&	&32	&17	&\multicolumn{3}{c |}{}\\
\hline
\end{tabular}
\end{center} 
\caption[]{Poisson $h$-multigrid with static condensation, two-grid results: iteration counts to achieve a residual tolerance of $10^{-10}$.  In each case, the fine mesh is identical to the the coarse mesh, once refined in $h$---in particular, the coarse and fine mesh have the same polynomial order.}
\label{table:PoissonHTwoGridResults}
\end{table}

The two-grid results for Stokes in two and three space dimensions can be found in Tables \ref{table:StokesPTwoGridResults} and \ref{table:StokesHTwoGridResults}.  As with Poisson, all the $p$-multigrid and $h$-multigrid operators scale, and the number of iterations required is relatively modest.

\begin{table}[h!b!p!]
\begin{center}
\begin{tabular}{| c  c  c  | c  c  c | c  c  c |}
\cline{1-9}
\multicolumn{9}{| c |}{$d=2$} \\
\cline{1-9}
$k$ & Mesh Width & Iterations	&$k$ & Mesh Width & Iterations	&$k$ & Mesh Width & Iterations	\\
1	&2	&16	&2	&2	&12	&4	&2	&14\\
	&4	&23	&	&4	&16	&	&4	&18\\
	&8	&24	&	&8	&17	&	&8	&19\\
	&16	&24	&	&16	&16	&	&16	&19\\
	&32	&24	&	&32	&16	&	&32	&19\\
	&64	&24	&	&64	&16	&	&64	&19\\
\hline
\multicolumn{9}{| c |}{$d=3$} \\
\hline 
$k$ & Mesh Width & Iterations	&$k$ & Mesh Width & Iterations	&\multicolumn{3}{c |}{}\\
1	&2	&34	&2	&2	&23	&\multicolumn{3}{c |}{}\\
	&4	&47	&	&4	&29	&\multicolumn{3}{c |}{}\\
	&8	&46	&	&8	&29	&\multicolumn{3}{c |}{}\\
	&16	&44	&	&16	&23	&\multicolumn{3}{c |}{}\\
	&32	&42	&	&32	&22	&\multicolumn{3}{c |}{}\\
\hline
\end{tabular}
\end{center} 
\caption[]{Stokes $p$-multigrid with static condensation, two-grid results: iteration counts to achieve a residual tolerance of $10^{-10}$.  In each case, the coarse mesh is geometrically identical to the fine, but with polynomial order half that of the fine mesh ($0$ when the fine mesh has polynomial order $1$).}
\label{table:StokesPTwoGridResults}
\end{table}

\begin{table}[h!b!p!]
\begin{center}
\begin{tabular}{| c  c  c  | c  c  c | c  c  c |}
\cline{1-9}
\multicolumn{9}{| c |}{$d=2$} \\
\cline{1-9}
$k$ & Mesh Width & Iterations	&$k$ & Mesh Width & Iterations	&$k$ & Mesh Width & Iterations \\
1	&4	&16	&2	&4	&15	&4	&4	&15\\
	&8	&18	&	&8	&17	&	&8	&17\\
	&16	&20	&	&16	&17	&	&16	&16\\
	&32	&20	&	&32	&16	&	&32	&16\\
	&64	&21	&	&64	&16	&	&64	&16\\
\hline
\multicolumn{9}{| c |}{$d=3$} \\
\hline 
$k$ & Mesh Width & Iterations		&$k$ & Mesh Width & Iterations		&\multicolumn{3}{c |}{}\\
1	&4	&27	&2	&4	&23	&\multicolumn{3}{c |}{}\\
	&8	&31	&	&8	&24	&\multicolumn{3}{c |}{}\\
	&16	&31	&	&16	&20	&\multicolumn{3}{c |}{}\\
	&32	&30	&\multicolumn{3}{| c |}{} &\multicolumn{3}{c |}{}\\
\hline
\end{tabular}
\end{center} 
\caption[]{Stokes $h$-multigrid with static condensation, two-grid results: iteration counts to achieve a residual tolerance of $10^{-10}$.  In each case, the fine mesh is identical to the the coarse mesh, once refined in $h$---in particular, the coarse and fine mesh have the same polynomial order.}
\label{table:StokesHTwoGridResults}
\end{table}

Navier-Stokes is of particular interest because it involves spatially varying material data.  After the first Newton step in Navier-Stokes, the background flow will be non-zero, and therefore the material data will vary in space.  To focus on this case, we start with a linear mesh and take three Newton steps; we use this as the background flow for the solve on our fine mesh---this is the solve for which we report iteration counts.

As is our default throughout this paper, in the first set of experiments we use a test space enrichment $\Delta k = d$, where $d=2$ is the spatial dimension.  The results for $p$ and $h$ operators can be found in Tables \ref{table:NavierStokesPTwoGridResultsDeltak2} and \ref{table:NavierStokesHTwoGridResultsDeltak2}.  Here, the multigrid preconditioner appears to be robust in $h$, but not in $p$ (though in the $p$ case the iteration counts do not grow by too much); the $k=4$ results have higher iteration counts than do the $k=1$ results, for example.

\begin{table}[h!b!p!]
\begin{center}
\begin{tabular}{| c  c  c  | c  c  c | c  c  c |}
\cline{1-9}
\multicolumn{9}{| c |}{$d=2$} \\
\cline{1-9}
$k$ & Mesh Width & Iterations	&$k$ & Mesh Width & Iterations	&$k$ & Mesh Width & Iterations	\\
1	&2	&18	&2	&2	&18	&4	&2	&21\\
	&4	&32	&	&4	&26	&	&4	&32\\
	&8	&37	&	&8	&29	&	&8	&35\\
	&16	&36	&	&16	&29	&	&16	&38\\
	&32	&33	&	&32	&29	&	&32	&39\\
	&64	&30	&	&64	&28	&	&64	&40\\
\hline
\end{tabular}
\end{center} 
\caption[]{Navier-Stokes $p$-multigrid with static condensation, two-grid results with $\Delta k=d=2$: iteration counts to achieve a residual tolerance of $10^{-10}$.  In each case, the coarse mesh is geometrically identical to the fine, but with polynomial order half that of the fine mesh ($0$ when the fine mesh has polynomial order $1$).}
\label{table:NavierStokesPTwoGridResultsDeltak2}
\end{table}

\begin{table}[h!b!p!]
\begin{center}
\begin{tabular}{| c  c  c  | c  c  c | c  c  c |}
\cline{1-9}
\multicolumn{9}{| c |}{$d=2$} \\
\cline{1-9}
$k$ & Mesh Width & Iterations	&$k$ & Mesh Width & Iterations	&$k$ & Mesh Width & Iterations	\\
1	&4	&19	&2	&4	&20	&4	&4	&21\\
	&8	&24	&	&8	&24	&	&8	&21\\
	&16	&23	&	&16	&22	&	&16	&21\\
	&32	&23	&	&32	&22	&	&32	&20\\
	&64	&24	&	&64	&24	&	&64	&23\\
\hline
\end{tabular}
\end{center} 
\caption[]{Navier-Stokes $h$-multigrid with static condensation, two-grid results with $\Delta k=d=2$: iteration counts to achieve a residual tolerance of $10^{-10}$.  In each case, the fine mesh is identical to the the coarse mesh, once refined in $h$---in particular, the coarse and fine mesh have the same polynomial order.}
\label{table:NavierStokesHTwoGridResultsDeltak2}
\end{table}

If we repeat our experiment, now with $\Delta k = 4$, we get essentially the same results for the $h$ preconditioners, but markedly better results for the $p$ preconditioners, as can be seen in Tables \ref{table:NavierStokesPTwoGridResultsDeltak4} and \ref{table:NavierStokesHTwoGridResultsDeltak4}.  Here, for both sets of preconditioners, we see results that closely resemble what we saw for Stokes and Poisson: the higher-order meshes have generally lower iteration counts, and the finer meshes at a given polynomial order have a roughly fixed iteration count.  

\clearpage
It is straightforward to show that DPG is equivalent to a non-standard mixed formulation involving the test space \cite{demkowicz2015ices}.    Dahmen et al.\ showed that the convergence of the Uzawa iteration for this mixed formulation depends on the approximation of the test space \cite{dahmen2012adaptive}.  This suggests that increasing the degree of enrichment $\Delta k$ for the DPG test functions improves the effectiveness of the preconditioner; this was observed independently by Gopalakrishnan\footnote{Private communication.} in the preconditioning of DPG for Maxwell's equations.

\begin{table}[h!b!p!]
\begin{center}
\begin{tabular}{| c  c  c  | c  c  c | c  c  c |}
\cline{1-9}
\multicolumn{9}{| c |}{$d=2$} \\
\cline{1-9}
$k$ & Mesh Width & Iterations	&$k$ & Mesh Width & Iterations	&$k$ & Mesh Width & Iterations	\\
1	&2	&17	&2	&2	&16	&4	&2	&17\\
	&4	&28	&	&4	&20	&	&4	&20\\
	&8	&33	&	&8	&20	&	&8	&21\\
	&16	&31	&	&16	&19	&	&16	&21\\
	&32	&32	&	&32	&20	&	&32	&21\\
	&64	&28	&	&64	&19	&	&64	&21\\
\hline
\end{tabular}
\end{center} 
\caption[]{Navier-Stokes $p$-multigrid with static condensation, two-grid results with $\Delta k=4$: iteration counts to achieve a residual tolerance of $10^{-10}$.  In each case, the coarse mesh is geometrically identical to the fine, but with polynomial order half that of the fine mesh ($0$ when the fine mesh has polynomial order $1$).}
\label{table:NavierStokesPTwoGridResultsDeltak4}
\end{table}

\begin{table}[h!b!p!]
\begin{center}
\begin{tabular}{| c  c  c  | c  c  c | c  c  c |}
\cline{1-9}
\multicolumn{9}{| c |}{$d=2$} \\
\cline{1-9}
$k$ & Mesh Width & Iterations	&$k$ & Mesh Width & Iterations	&$k$ & Mesh Width & Iterations	\\
1	&4	&18	&2	&4	&18	&4	&4	&18\\
	&8	&22	&	&8	&19	&	&8	&20\\
	&16	&22	&	&16	&19	&	&16	&19\\
	&32	&23	&	&32	&20	&	&32	&19\\
	&64	&24	&	&64	&22	&	&64	&19\\
\hline
\end{tabular}
\end{center} 
\caption[]{Navier-Stokes $h$-multigrid with static condensation, two-grid results with $\Delta k=4$: iteration counts to achieve a residual tolerance of $10^{-10}$.  In each case, the fine mesh is identical to the the coarse mesh, once refined in $h$---in particular, the coarse and fine mesh have the same polynomial order.}
\label{table:NavierStokesHTwoGridResultsDeltak4}
\end{table}

\subsection{Multigrid Experiments}
\paragraph{Uniform Refinements}
For the uniform refinement cases, we use the same problems as described in Section \ref{sec:problemDefinitions} for Stokes and Navier-Stokes, and construct mesh hierarchies as described in Section \ref{sec:meshHierarchy}.  We are particularly interested in the effect of the number of grids in the hierarchy on the iteration count.

Results for Stokes with the full $p$-hierarchy (i.e., with parameter \texttt{skipIntermediateP} is taken to be false) can be found in Table \ref{table:StokesMultigridResults}; the iteration counts do grow as the number of grids increases, with more rapid growth as the number of $h$ levels increases, but the total iteration counts remain modest.  Results with \texttt{skipIntermediateP = true} are shown in Table \ref{table:StokesMultigridResultsSkipIntermediate}.  In most cases, the iteration counts are slightly higher with this option.  However, this is the computationally cheaper case: fewer smoother and prolongation operators need to be computed and stored; this is the approach we advocate in practice.

As suggested by our results above, we run Navier-Stokes both for $\Delta k = 2$ and for $\Delta k = k_{\rm fine}$.  We again begin by taking \texttt{skipIntermediateP} to be false. The results for $\Delta k = 2$ can be found in Tables \ref{table:NavierStokesMultigridResultsDeltaK2}; here, the iteration counts grow more rapidly in the number of multigrid levels than they do for Stokes.  Table \ref{table:NavierStokesMultigridResultsDeltaKFine} shows the results for $\Delta k = k_{\rm fine}$; here, the higher-order results are roughly in line with those for Stokes.

Repeating the Navier-Stokes experiments with \texttt{skipIntermediateP=true}, we find that iteration counts are again somewhat higher than for the alternative---the results are shown in Tables \ref{table:NavierStokesMultigridResultsDeltaK2SkipIntermediate} and \ref{table:NavierStokesMultigridResultsDeltaKFineSkipIntermediate}.  Again we see a significant reduction in iteration count by taking $\Delta k = k_{\rm fine}$.

\begin{table}[h!b!p!]
\begin{center}
\begin{tabular}{| c  c  c  c  c  c c |}
\hline
$k_{\rm fine}$ & $k_{\rm coarse}$ & $k$ levels & Fine Mesh Width & Coarse Mesh Width & $h$ levels & Iterations\\
                   1                   &1                   &0                   &4                   &2                   &1                   &19\\
                   1                   &1                   &0                   &8                   &2                   &2                   &28\\
                   1                   &1                   &0                   &16                   &2                   &3                   &36\\
                   1                   &1                   &0                   &32                   &2                   &4                   &48\\
\hline
                   2                   &1                   &1                   &2                   &2                   &0                   &15\\
                   4                   &1                   &2                   &2                   &2                   &0                   &21\\
                   8                   &1                   &3                   &2                   &2                   &0                   &26\\
                  16                   &1                   &4                   &2                   &2                   &0                   &31\\
\hline
                   4                   &1                   &2                   &4                   &2                   &1                   &33\\
                   4                   &1                   &2                   &8                   &2                   &2                   &39\\
                   4                   &1                   &2                   &16                   &2                   &3                   &44\\
                   4                   &1                   &2                   &32                   &2                   &4                   &54\\
\hline
                   2                   &1                   &1                   &32                   &32                   &0                   &18\\
                   4                   &1                   &2                   &32                   &32                   &0                   &25\\
\hline
\end{tabular}
\end{center} 
\caption[]{Stokes solves for multigrid (2D).  We see some growth in iteration counts as the number of $h$ and $k$ levels increases; the growth is more pronounced for the $h$-level increases.}
\label{table:StokesMultigridResults}
\end{table}

\begin{table}[h!b!p!]
\begin{center}
\begin{tabular}{| c  c  c  c  c  c c |}
\hline
$k_{\rm fine}$ & $k_{\rm coarse}$ & $k$ levels & Fine Mesh Width & Coarse Mesh Width & $h$ levels & Iterations\\
                   4                   &1                   &1                   &2                   &2                   &0                   &22\\
                   8                   &1                   &1                   &2                   &2                   &0                   &30\\
                  16                   &1                   &1                   &2                   &2                   &0                   &39\\
\hline
                   4                   &1                   &1                   &4                   &2                   &1                   &36\\
                   4                   &1                   &1                   &8                   &2                   &2                   &41\\
                   4                   &1                   &1                   &16                   &2                   &3                   &46\\
                   4                   &1                   &1                   &32                   &2                   &4                   &54\\
\hline
                   4                   &1                   &1                   &32                   &32                   &0                   &27\\
\hline
\end{tabular}
\end{center} 
\caption[]{Stokes solves for multigrid (2D), skipping intermediate polynomial orders between $k_{\rm coarse}$ and $k_{\rm fine}$.}
\label{table:StokesMultigridResultsSkipIntermediate}
\end{table}

\begin{table}[h!b!p!]
\begin{center}
\begin{tabular}{| c  c  c  c  c  c c |}
\hline
$k_{\rm fine}$ & $k_{\rm coarse}$ & $k$ levels & Fine Mesh Width & Coarse Mesh Width & $h$ levels & Iterations\\
                   1                   &1                   &0                   &4                   &2                   &1                   &20\\
                   1                   &1                   &0                   &8                   &2                   &2                   &28\\
                   1                   &1                   &0                   &16                   &2                   &3                   &32\\
                   1                   &1                   &0                   &32                   &2                   &4                   &43\\
\hline
                   2                   &1                   &1                   &2                   &2                   &0                   &23\\
                   4                   &1                   &2                   &2                   &2                   &0                   &33\\
                   8                   &1                   &3                   &2                   &2                   &0                   &46\\
                  16                   &1                   &4                   &2                   &2                   &0                   &59\\
\hline
                   4                   &1                   &2                   &4                   &2                   &1                   &54\\
                   4                   &1                   &2                   &8                   &2                   &2                   &63\\
                   4                   &1                   &2                   &16                   &2                   &3                   &80\\
                   4                   &1                   &2                   &32                   &2                   &4                   &94\\
\hline
                   2                   &1                   &1                   &32                   &32                   &0                   &33\\
                   4                   &1                   &2                   &32                   &32                   &0                   &62\\
\hline
\end{tabular}
\end{center} 
\caption[]{Navier-Stokes solves for multigrid (2D), with $\Delta k = 2$.  The growth in iteration counts is more pronounced than it is for the Stokes results reported in Table \ref{table:StokesMultigridResults}.}
\label{table:NavierStokesMultigridResultsDeltaK2}
\end{table}

\begin{table}[h!b!p!]
\begin{center}
\begin{tabular}{| c  c  c  c  c  c c |}
\hline
$k_{\rm fine}$ & $k_{\rm coarse}$ & $k$ levels & Fine Mesh Width & Coarse Mesh Width & $h$ levels & Iterations\\
                   1                   &1                   &0                   &4                   &2                   &1                   &25\\
                   1                   &1                   &0                   &8                   &2                   &2                   &42\\
                   1                   &1                   &0                   &16                   &2                   &3                   &50\\
                   1                   &1                   &0                   &32                   &2                   &4                   &68\\
\hline
                   2                   &1                   &1                   &2                   &2                   &0                   &23\\
                   4                   &1                   &2                   &2                   &2                   &0                   &26\\
                   8                   &1                   &3                   &2                   &2                   &0                   &22\\
\hline
                   4                   &1                   &2                   &4                   &2                   &1                   &34\\
                   4                   &1                   &2                   &8                   &2                   &2                   &39\\
                   4                   &1                   &2                   &16                   &2                   &3                   &45\\
                   4                   &1                   &2                   &32                   &2                   &4                   &52\\
\hline
                   2                   &1                   &1                   &32                   &32                   &0                   &33\\
                   4                   &1                   &2                   &32                   &32                   &0                   &29\\
\hline
\end{tabular}
\end{center} 
\caption[]{Navier-Stokes solves for multigrid (2D), with $\Delta k = k_{\rm fine}$.  For the higher-order meshes, the growth in iteration counts is considerably reduced compared to the $\Delta k = 2$ results reported in Table \ref{table:NavierStokesMultigridResultsDeltaK2}, and roughly in line with those of the Stokes results in Table \ref{table:StokesMultigridResults}.}
\label{table:NavierStokesMultigridResultsDeltaKFine}
\end{table}

\begin{table}[h!b!p!]
\begin{center}
\begin{tabular}{| c  c  c  c  c  c c |}
\hline
$k_{\rm fine}$ & $k_{\rm coarse}$ & $k$ levels & Fine Mesh Width & Coarse Mesh Width & $h$ levels & Iterations\\
                   4                   &1                   &1                   &2                   &2                   &0                   &36\\
                   8                   &1                   &1                   &2                   &2                   &0                   &48\\
                  16                   &1                   &1                   &2                   &2                   &0                   &68\\
\hline
                   4                   &1                   &1                   &4                   &2                   &1                   &59\\
                   4                   &1                   &1                   &8                   &2                   &2                   &69\\
                   4                   &1                   &1                   &16                   &2                   &3                   &89\\
                   4                   &1                   &1                   &32                   &2                   &4                   &106\\
\hline
                   4                   &1                   &1                   &32                   &32                   &0                   &67\\
\hline
\end{tabular}
\end{center} 
\caption[]{Navier-Stokes solves for multigrid (2D), with $\Delta k = 2$, skipping intermediate polynomial orders between $k_{\rm coarse}$ and $k_{\rm fine}$.  Iteration counts are higher than for those in which we did not skip polynomial orders, shown in Table \ref{table:NavierStokesMultigridResultsDeltaK2}; however, the operators employed here are cheaper to compute and store than those in which we do not skip intermediate polynomial orders.}
\label{table:NavierStokesMultigridResultsDeltaK2SkipIntermediate}
\end{table}

\begin{table}[h!b!p!]
\begin{center}
\begin{tabular}{| c  c  c  c  c  c c |}
\hline
$k_{\rm fine}$ & $k_{\rm coarse}$ & $k$ levels & Fine Mesh Width & Coarse Mesh Width & $h$ levels & Iterations\\
                   4                   &1                   &1                   &2                   &2                   &0                   &27\\
                   8                   &1                   &1                   &2                   &2                   &0                   &25\\
\hline
                   4                   &1                   &1                   &4                   &2                   &1                   &37\\
                   4                   &1                   &1                   &8                   &2                   &2                   &43\\
                   4                   &1                   &1                   &16                   &2                   &3                   &48\\
                   4                   &1                   &1                   &32                   &2                   &4                   &57\\
\hline
                   4                   &1                   &1                   &32                   &32                   &0                   &33\\
\hline
\end{tabular}
\end{center} 
\caption[]{Navier-Stokes solves for multigrid (2D), with $\Delta k = k_{\rm fine}$, skipping intermediate polynomial orders between $k_{\rm coarse}$ and $k_{\rm fine}$.  Iteration counts are somewhat higher than for those in which we did not skip polynomial orders, shown in Table \ref{table:NavierStokesMultigridResultsDeltaKFine}; however, these operators are cheaper to compute and store than those.}
\label{table:NavierStokesMultigridResultsDeltaKFineSkipIntermediate}
\end{table}

\paragraph{Adaptive Refinements}
\label{sec:adaptiveCavityFlow}
An important motivation for using DPG is that one can obtain meaningful solutions even on a coarse mesh, and use the method's built-in error estimation to drive adaptive refinements.  Therefore, it is worth confirming that our solver behaves reasonably in the presence of hanging nodes.  Here, we consider Stokes and Navier-Stokes for the lid-driven cavity flow problem.

\begin{figure}[h!b!p!]
\begin{centering}
\includegraphics[width=.50\linewidth]{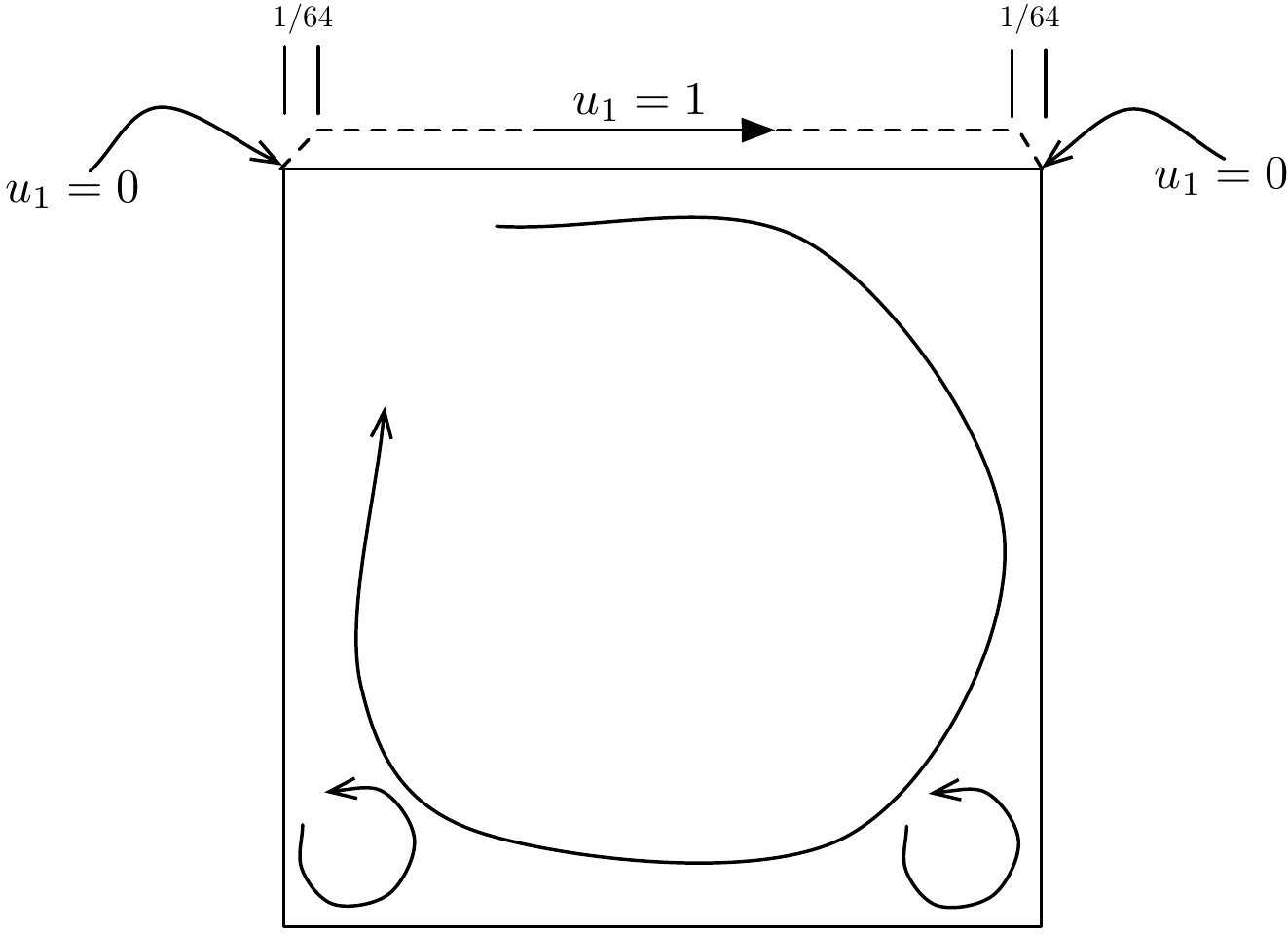}\\
\caption{Lid-driven cavity flow schematic.  The boundary conditions interpolate between $u_1=1$ at the lid and $u_1=0$ along the side walls.}
\label{fig:cavityFlowCartoon}
\end{centering}
\end{figure}

The cavity is defined on a domain $[0,1]^2$, with no-slip boundary conditions on the $y=0$, $x=0$, and $x=1$ walls; a schematic can be seen in Figure \ref{fig:cavityFlowCartoon}.  The lid at $y=1$ moves at a constant unit velocity.  Because physically there will be some continuous transition between the zero velocity at the wall and the unit velocity at the lid---and because a discontinuity in the velocity would imply an exact solution outside $H^1$---in our studies we employ a thin linear interpolation along the left and right sides of the lid.  The width of the transition we use is $\epsilon = \frac{1}{64}$.  A consequence of this choice is that we do not exactly capture the boundary conditions until the sixth refinement, when starting from a $2 \times 2$ initial mesh.

We use static condensation at every grid level.  This has only a minor effect on the iteration counts (often changing them not at all, and sometimes reducing or increasing them by 1 or 2), but results in a substantial reduction in runtime and memory requirements.  We use as the coarsest mesh a $2 \times 2$ mesh of polynomial order 1, and use the procedure discussed in Section \ref{sec:gmgDetails} to define the sequence of grids.  We perform CG iterations until a tolerance of $10^{-6}$ is met.  On each refinement step, we use as the initial guess the solution from the prior step.  (For comparison, we also show results using a zero initial guess.)

Refinements are performed using a greedy algorithm: at each refinement step, we compute the maximum element error $err_{\rm max}$, and refine all elements with error greater than 20\% of $err_{\rm max}$.  

\begin{figure}[h!b!p!]
\begin{centering}
\includegraphics[width=.90\linewidth]{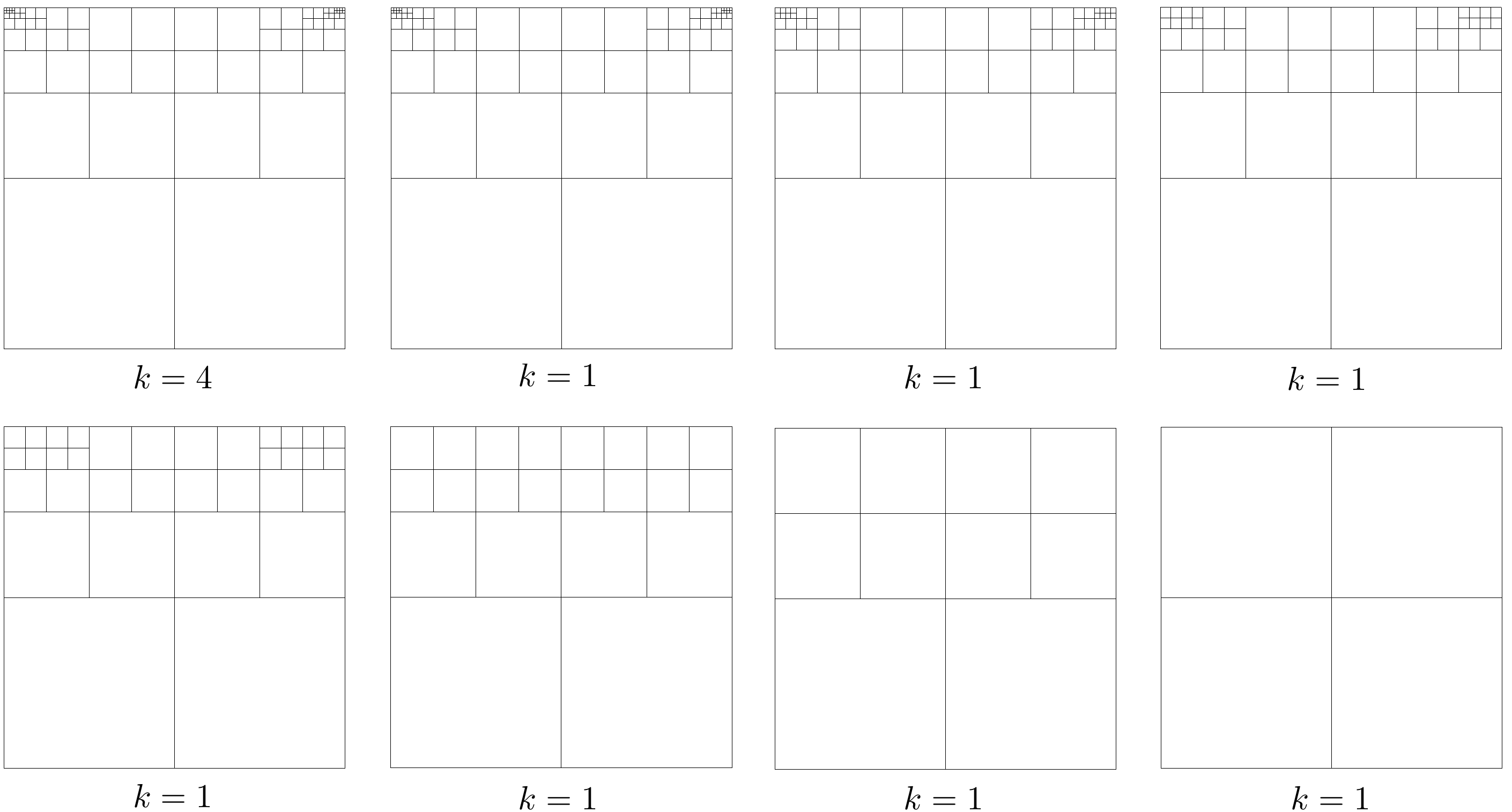}\\
\caption{Multigrid hierarchy for quartic lid-driven cavity flow with $k=4$ after six adaptive mesh refinements.  Here, the option \texttt{skipIntermediateP} described in Section \ref{sec:meshHierarchy} is taken to be true; if it were false, there would be an additional $k=2$ mesh between the $k=4$ mesh and the finest $k=1$ mesh.}
\label{fig:cavityFlowHierarchy}
\end{centering}
\end{figure}

We solve using mesh hierarchies defined according to the algorithm defined in Section \ref{sec:meshHierarchy}, with the option to skip intermediate polynomial orders set to \texttt{true}.  An example mesh hierarchy, for the sixth refinement of a quartic mesh, is shown in Figure \ref{fig:cavityFlowHierarchy}. The results at polynomial orders of $k=1,2,4,$ and $8$ are listed in Table \ref{table:StokesCavityAdaptiveIterationCounts}.  A few phenomena are worth noting.  First, the iteration counts required are considerably lower for higher-order meshes.  Second, using the solution from the previous mesh as an initial guess provides considerable benefit on higher-order meshes.  Finally, though in some cases the iteration counts grow as we refine the mesh, they do so relatively modestly, especially on higher-order meshes when using the previous solution as initial guess.

\begin{table}[h!b!p!]
\begin{center}
\begin{tabular}{ c  c  c  c  c  c  c  c  c }
\hline
\multicolumn{8}{ c }{$k=1$} \\
\hline 
Ref. \# & $h_{\rm max}$ & $h_{\rm min}$ & $\frac{h_{\rm max}}{h_{\rm min}}$ & Elements & Energy Error & Iter. &w/Zero Guess\\
                   0                   &1/2                   &1/2                   &1                   &4                   &8.96e-01                   &1                   &1\\
                   1                   &1/4                   &1/2                   &2                   &10                   &9.36e-01                   &8                   &9\\
                   2                   &1/8                   &1/2                   &4                   &16                   &9.23e-01                   &12                   &12\\
                   3                   &1/16                   &1/2                   &8                   &22                   &9.25e-01                   &16                   &16\\
                   4                   &1/32                   &1/2                   &16                   &28                   &7.72e-01                   &17                   &17\\
                   5                   &1/64                   &1/4                   &16                   &88                   &4.40e-01                   &34                   &35\\
                   6                   &1/128                   &1/4                   &32                   &106                   &2.67e-01                   &37                   &38\\
                   7                   &1/256                   &1/4                   &64                   &268                   &1.26e-01                   &56                   &60\\
                   8                   &1/512                   &1/8                   &64                   &358                   &8.38e-02                   &49                   &71\\
\hline
\multicolumn{8}{ c }{$k=2$} \\
\hline 
Ref. \# & $h_{\rm max}$ & $h_{\rm min}$ & $\frac{h_{\rm max}}{h_{\rm min}}$ & Elements & Energy Error & Iter. &w/Zero Guess\\
                   0                   &1/2                   &1/2                   &1                   &4                   &8.48e-01                   &11                   &11\\
                   1                   &1/4                   &1/2                   &2                   &10                   &8.12e-01                   &15                   &16\\
                   2                   &1/8                   &1/2                   &4                   &16                   &7.53e-01                   &18                   &19\\
                   3                   &1/16                   &1/2                   &8                   &22                   &6.64e-01                   &23                   &25\\
                   4                   &1/32                   &1/2                   &16                   &28                   &3.49e-01                   &16                   &28\\
                   5                   &1/64                   &1/2                   &32                   &34                   &2.02e-01                   &19                   &29\\
                   6                   &1/128                   &1/4                   &32                   &100                   &8.72e-02                   &32                   &49\\
                   7                   &1/256                   &1/4                   &64                   &124                   &4.97e-02                   &25                   &54\\
                   8                   &1/512                   &1/4                   &128                   &178                   &2.85e-02                   &25                   &60\\
\hline
\multicolumn{8}{ c }{$k=4$} \\
\hline
Ref. \# & $h_{\rm max}$ & $h_{\rm min}$ & $\frac{h_{\rm max}}{h_{\rm min}}$ & Elements & Energy Error & Iter. &w/Zero Guess\\
                   0                   &1/2                   &1/2                   &1                   &4                   &7.27e-01                   &15                   &15\\
                   1                   &1/4                   &1/2                   &2                   &10                   &6.30e-01                   &20                   &22\\
                   2                   &1/8                   &1/2                   &4                   &16                   &5.82e-01                   &15                   &28\\
                   3                   &1/16                   &1/2                   &8                   &22                   &2.06e-01                   &20                   &41\\
                   4                   &1/32                   &1/2                   &16                   &28                   &7.54e-02                   &21                   &46\\
                   5                   &1/64                   &1/2                   &32                   &34                   &5.90e-02                   &18                   &33\\
                   6                   &1/128                   &1/2                   &64                   &70                   &3.01e-02                   &21                   &55\\
                   7                   &1/256                   &1/2                   &128                   &88                   &1.55e-02                   &19                   &63\\
                   8                   &1/512                   &1/2                   &256                   &106                   &8.63e-03                   &15                   &70\\
\hline
\multicolumn{8}{ c }{$k=8$} \\
\hline
Ref. \# & $h_{\rm max}$ & $h_{\rm min}$ & $\frac{h_{\rm max}}{h_{\rm min}}$ & Elements & Energy Error & Iter. &w/Zero Guess\\
                   0                   &1/2                   &1/2                   &1                   &4                   &4.72e-01                   &20                   &20\\
                   1                   &1/4                   &1/2                   &2                   &10                   &3.53e-01                   &16                   &29\\
                   2                   &1/8                   &1/2                   &4                   &16                   &1.29e-01                   &14                   &38\\
                   3                   &1/16                   &1/2                   &8                   &22                   &9.97e-02                   &19                   &54\\
                   4                   &1/32                   &1/2                   &16                   &28                   &2.98e-02                   &19                   &60\\
                   5                   &1/64                   &1/2                   &32                   &34                   &1.79e-02                   &19                   &38\\
                   6                   &1/128                   &1/2                   &64                   &70                   &9.06e-03                   &12                   &42\\
                   7                   &1/256                   &1/2                   &128                   &91                   &4.58e-03                   &10                   &52\\
                   8                   &1/512                   &1/2                   &256                   &112                   &2.27e-03                   &9                   &58\\
\hline 
\end{tabular}
\end{center} 
\caption[]{Stokes cavity flow: iteration counts to achieve a residual tolerance of $10^{-6}$.  The rightmost column lists the iteration counts with a zero initial guess on each refined mesh, while the column immediately left of that lists iteration counts starting from an initial guess corresponding to the solution on the mesh of the previous refinement step.}
\label{table:StokesCavityAdaptiveIterationCounts}
\end{table}

For Navier-Stokes, following the approach described in \cite{DPGNavierStokes}, we define an initial nonlinear stopping threshold $\epsilon_0 = 10^{-4}$ for the initial mesh; we perform Newton iterations until the $L^2$ norm of the field variables in the solution increment is less than $\epsilon_0$.  We define the relative error $e^i_{\rm rel}$ for solution $u^i_h$ in terms of the DPG energy error $\norm{e^i}_E$ and the energy of the solution $\norm{u^{i}_h}_E$:
\begin{align*}
e^i_{\rm rel} = \frac{\norm{e^i}_E}{\norm{u^i_h}_E} = \frac{\norm{u - u^i_h}_E}{\norm{u^i_h}_E} = \frac{\norm{l(\cdot) - b(u^{i}_h,\cdot)}_{V^{'}_h}}{\norm{b(u^i_h,\cdot)}_{V^{'}_h}}.
\end{align*}
Note that the bilinear form $b(\cdot,\cdot)$ is for the linearized problem, and therefore depends on the background flow.  Because the graph norm we employ on the test space depends on the bilinear form, the norms $\norm{\cdot}_E$ and $\norm{\cdot}_{V^{'}_h}$ also depend on the background flow.

Results for ${\rm Re} = 100$ with $k=1,2,4$ with $\Delta k=2$, a conjugate gradient stopping tolerance of $10^{-6}$, and \texttt{skipIntermediateP=true} are shown in Table \ref{table:NavierStokesCavityAdaptiveIterationCounts}.  The iteration counts for each Newton step are on average a little more than what we saw with Stokes when using a zero initial guess for the conjugate gradient iteration on each refinement step, but not dramatically so.  For the $k=8$ case, it appears that using $\Delta k=2$ does not suffice to resolve the optimal test functions; the evidence for this is that the energy norm increases substantially under refinement on finer meshes.  When we instead use $\Delta k=4$ and a conjugate gradient stopping tolerance of $10^{-9}$ in the $k=8$ case, we get the results shown in Table \ref{table:NavierStokesCavityAdaptiveIterationCountsK8}.

Results for ${\rm Re} = 1000$ with $k=1,2,4$ and $\Delta k=4$ are shown in Table \ref{table:NavierStokesCavityAdaptiveIterationCountsRe1e3}.  For finer meshes, the iteration counts per Newton step approach 300; this illustrates the limitations of the present black-box approach.

\begin{table}[h!b!p!]
\begin{center}
\begin{tabular}{ c  c  c  c  c  c  c  c  c}
\hline
\multicolumn{9}{ c }{$k=1$} \\
\hline
             Ref. \#     & $h_{\rm max}$     & $h_{\rm min}$& $\frac{h_{\rm max}}{h_{\rm min}}$          & Elements       & Energy Err.   & Nonlinear Steps  & Total Iterations       & Per Step\\
                   0                   &1/2                   &1/2                   &1                   &4                   &5.54e-01                   &7                   &7                   &1\\
                   1                   &1/4                   &1/4                   &1                   &16                   &4.65e-01                   &6                   &81                   &14\\
                   2                   &1/8                   &1/4                   &2                   &31                   &4.84e-01                   &5                   &99                   &20\\
                   3                   &1/16                   &1/4                   &4                   &46                   &6.72e-01                   &4                   &88                   &22\\
                   4                   &1/32                   &1/4                   &8                   &67                   &7.32e-01                   &5                   &111                   &22\\
                   5                   &1/64                   &1/4                   &16                   &94                   &5.73e-01                   &4                   &89                   &22\\
                   6                   &1/128                   &1/4                   &32                   &121                   &4.67e-01                   &3                   &70                   &23\\
                   7                   &1/256                   &1/4                   &64                   &274                   &2.86e-01                   &5                   &174                   &35\\
                   8                   &1/512                   &1/4                   &128                   &427                   &1.76e-01                   &4                   &178                   &45\\
\hline
\multicolumn{9}{ c }{$k=2$} \\
\hline
             Ref. \#     & $h_{\rm max}$     & $h_{\rm min}$& $\frac{h_{\rm max}}{h_{\rm min}}$          & Elements       & Energy Err.   & Nonlinear Steps  & Total Iterations       & Per Step\\
                   0                   &1/2                   &1/2                   &1                   &4                   &3.50e-01                   &7                   &99                   &14\\
                   1                   &1/4                   &1/2                   &2                   &10                   &2.44e-01                   &5                   &123                   &25\\
                   2                   &1/8                   &1/2                   &4                   &25                   &2.21e-01                   &5                   &217                   &43\\
                   3                   &1/16                   &1/2                   &8                   &34                   &2.37e-01                   &4                   &184                   &46\\
                   4                   &1/32                   &1/2                   &16                   &55                   &1.27e-01                   &5                   &256                   &51\\
                   5                   &1/64                   &1/4                   &16                   &103                   &1.42e-01                   &4                   &266                   &67\\
                   6                   &1/128                   &1/4                   &32                   &130                   &7.50e-02                   &4                   &248                   &62\\
                   7                   &1/256                   &1/4                   &64                   &247                   &3.95e-02                   &4                   &331                   &83\\
                   8                   &1/512                   &1/4                   &128                   &385                   &2.13e-02                   &4                   &331                   &83\\
\hline
\multicolumn{9}{ c }{$k=4$} \\
\hline
             Ref. \#     & $h_{\rm max}$     & $h_{\rm min}$& $\frac{h_{\rm max}}{h_{\rm min}}$          & Elements       & Energy Err.   & Nonlinear Steps  & Total Iterations       & Per Step\\
                   0                   &1/2                   &1/2                   &1                   &4                   &1.64e-01                   &6                   &145                   &24\\
                   1                   &1/4                   &1/2                   &2                   &10                   &1.29e-01                   &5                   &206                   &41\\
                   2                   &1/8                   &1/2                   &4                   &16                   &1.19e-01                   &5                   &270                   &54\\
                   3                   &1/16                   &1/2                   &8                   &28                   &2.83e-02                   &4                   &356                   &89\\
                   4                   &1/32                   &1/2                   &16                   &55                   &2.11e-02                   &5                   &541                   &108\\
                   5                   &1/64                   &1/2                   &32                   &79                   &2.13e-02                   &4                   &498                   &125\\
                   6                   &1/128                   &1/4                   &32                   &112                   &9.98e-03                   &4                   &544                   &136\\
                   7                   &1/256                   &1/4                   &64                   &160                   &4.98e-03                   &4                   &586                   &147\\
                   8                   &1/512                   &1/4                   &128                   &202                   &2.67e-03                   &4                   &511                   &128\\
\hline
\end{tabular}
\end{center} 
\caption[]{Navier-Stokes cavity flow for ${\rm Re} = 100$ with $\Delta k=2$: iteration counts to achieve a residual tolerance of $10^{-6}$.  We take Newton steps until the $L^2$ norm of the solution increment is below an adaptive tolerance which starts at $10^{-4}$, and grows tighter as the relative energy error of the previous refinement diminishes.  The reported energy errors are relative to the energy norm of the background flow.}
\label{table:NavierStokesCavityAdaptiveIterationCounts}
\end{table}

\begin{table}[h!b!p!]
\begin{center}
\begin{tabular}{ c  c  c  c  c  c  c  c  c}
\hline
\multicolumn{9}{ c }{$k=8$} \\
\hline
             Ref. \#     & $h_{\rm max}$     & $h_{\rm min}$& $\frac{h_{\rm max}}{h_{\rm min}}$          & Elements       & Energy Err.   & Nonlinear Steps  & Total Iterations       & Per Step\\
                   0                   &1/2                   &1/2                   &1                   &4                   &4.08e-01                   &5                   &136                   &27\\
                   1                   &1/4                   &1/2                   &2                   &10                   &2.68e-01                   &5                   &212                   &42\\
                   2                   &1/8                   &1/2                   &4                   &16                   &8.74e-02                   &4                   &239                   &60\\
                   3                   &1/16                   &1/2                   &8                   &34                   &6.69e-02                   &4                   &402                   &101\\
                   4                   &1/32                   &1/2                   &16                   &46                   &2.29e-02                   &4                   &446                   &112\\
                   5                   &1/64                   &1/2                   &32                   &58                   &3.10e-02                   &5                   &580                   &116\\
                   6                   &1/128                   &1/2                   &64                   &88                   &1.55e-02                   &4                   &548                   &137\\
                   7                   &1/256                   &1/2                   &128                   &106                   &9.68e-03                   &3                   &419                   &140\\
                   8                   &1/512                   &1/4                   &128                   &241                   &3.66e-03                   &3                   &591                   &197\\
\hline
\end{tabular}
\end{center} 
\caption[]{Navier-Stokes cavity flow for ${\rm Re} = 100$ with $k=8$ and $\Delta k=4$: iteration counts to achieve a residual tolerance of $10^{-9}$.  We take Newton steps until the $L^2$ norm of the solution increment is below an adaptive tolerance which starts at $10^{-4}$, and grows tighter as the relative energy error of the previous refinement diminishes.  The reported energy errors are relative to the energy norm of the background flow.}
\label{table:NavierStokesCavityAdaptiveIterationCountsK8}
\end{table}

\begin{table}[h!b!p!]
\begin{center}
\begin{tabular}{ c  c  c  c  c  c  c  c  c}
\hline
\multicolumn{9}{ c }{$k=1$} \\
\hline
           Ref. \#     & $h_{\rm max}$     & $h_{\rm min}$& $\frac{h_{\rm max}}{h_{\rm min}}$          & Elements       & Energy Err.   & Nonlinear Steps  & Total Iterations       & Per Step\\
                   0                   &1/2                   &1/2                   &1                   &4                   &6.42e-01                   &4                   &4                   &1\\
                   1                   &1/4                   &1/4                   &1                   &16                   &4.67e-01                   &9                   &166                   &18\\
                   2                   &1/8                   &1/4                   &2                   &43                   &5.19e-01                   &8                   &342                   &43\\
                   3                   &1/16                   &1/4                   &4                   &85                   &6.69e-01                   &30                   &1884                   &63\\
                   4                   &1/32                   &1/4                   &8                   &142                   &8.04e-01                   &30                   &2499                   &83\\
                   5                   &1/64                   &1/4                   &16                   &172                   &1.01e+00                   &6                   &466                   &78\\
                   6                   &1/128                   &1/4                   &32                   &211                   &7.75e-01                   &3                   &237                   &79\\
                   7                   &1/256                   &1/4                   &64                   &262                   &5.31e-01                   &3                   &237                   &79\\
                   8                   &1/512                   &1/4                   &128                   &418                   &3.74e-01                   &3                   &286                   &95\\
\hline
\multicolumn{9}{ c }{$k=2$} \\
\hline
             Ref. \#     & $h_{\rm max}$     & $h_{\rm min}$& $\frac{h_{\rm max}}{h_{\rm min}}$          & Elements       & Energy Err.   & Nonlinear Steps  & Total Iterations       & Per Step\\
                   0                   &1/2                   &1/2                   &1                   &4                   &4.30e-01                   &7                   &129                   &18\\
                   1                   &1/4                   &1/2                   &2                   &13                   &4.84e-01                   &30                   &1438                   &48\\
                   2                   &1/8                   &1/2                   &4                   &31                   &2.52e-01                   &12                   &1125                   &94\\
                   3                   &1/16                   &1/2                   &8                   &61                   &1.73e-01                   &5                   &694                   &139\\
                   4                   &1/32                   &1/4                   &8                   &151                   &1.10e-01                   &10                   &1873                   &187\\
                   5                   &1/64                   &1/4                   &16                   &253                   &1.73e-01                   &11                   &2847                   &259\\
                   6                   &1/128                   &1/4                   &32                   &265                   &1.56e-01                   &21                   &5380                   &256\\
                   7                   &1/256                   &1/4                   &64                   &346                   &1.20e-01                   &21                   &5701                   &271\\
                   8                   &1/512                   &1/4                   &128                   &409                   &8.57e-02                   &30                   &8594                   &286\\
\hline
\multicolumn{9}{ c }{$k=4$} \\
\hline
             Ref. \#     & $h_{\rm max}$     & $h_{\rm min}$& $\frac{h_{\rm max}}{h_{\rm min}}$          & Elements       & Energy Err.   & Nonlinear Steps  & Total Iterations       & Per Step\\
                   0                   &1/2                   &1/2                   &1                   &4                   &2.54e-01                   &20                   &797                   &40\\
                   1                   &1/4                   &1/2                   &2                   &10                   &1.29e-01                   &7                   &547                   &78\\
                   2                   &1/8                   &1/2                   &4                   &25                   &7.35e-02                   &6                   &847                   &141\\
                   3                   &1/16                   &1/4                   &4                   &67                   &3.49e-02                   &7                   &1677                   &240\\
                   4                   &1/32                   &1/4                   &8                   &88                   &5.42e-02                   &4                   &1070                   &268\\
                   5                   &1/64                   &1/4                   &16                   &109                   &2.08e-02                   &4                   &1093                   &273\\
                   6                   &1/128                   &1/4                   &32                   &142                   &3.46e-02                   &4                   &1136                   &284\\
                   7                   &1/256                   &1/4                   &64                   &175                   &2.15e-02                   &5                   &1476                   &295\\
                   8                   &1/512                   &1/4                   &128                   &208                   &1.60e-02                   &4                   &1179                   &295\\
\hline
\end{tabular}
\end{center} 
\caption[]{Navier-Stokes cavity flow for ${\rm Re} = 1000$ with $k=1,2,4$ and $\Delta k=4$: iteration counts to achieve a residual tolerance of $10^{-9}$.  We take up to 30 Newton steps until the $L^2$ norm of the solution increment is below an adaptive tolerance which starts at $10^{-2}$, and grows tighter as the relative energy error of the previous refinement diminishes.  The reported energy errors are relative to the energy norm of the background flow.}
\label{table:NavierStokesCavityAdaptiveIterationCountsRe1e3}
\end{table}

\section{Some Timing Results: Scaling of the Camellia Implementation}
\label{sec:scalability}
As discussed in the introduction, a primary driver for our efforts is development of DPG solvers that scale to large machines.  With that in mind, in this section we examine timing results for the same Stokes problem with smooth solution as we examined in Section \ref{sec:NumericalExperiments}, running on Argonne's Mira supercomputer.

We ran on between 512 and 4096 nodes on Mira, using 8 MPI ranks per node.  This afforded us 2 GB of RAM per MPI rank,\footnote{We required approximately 1.7 GB of RAM per rank for the smallest node count, where we had 8 fine elements per rank; for the largest node count (with one fine element per rank), we required about 750 MB per rank.  One significant component of the memory cost is the storage of the factored Schwarz blocks for the smoother at each grid level; the cost of these is approximately the same as storing the system matrix.  When using static condensation, as we do here, the memory cost of both the Schwarz blocks and the system matrix are reduced, but to reduce compute time we store the local (uncondensed) stiffness matrices, and these take as much memory as the global uncondensed stiffness matrix would.} and two Blue Gene/Q cores per rank.  In each case, we solved on a 32,768-element quartic mesh, with a 512-element coarse mesh of lowest (constant) order.  The fine mesh has 76 million degrees of freedom total, with 14 million trace degrees of freedom; the coarsest mesh has 25,000 total degrees of freedom, with 7,400 trace degrees of freedom.  As we have done throughout, here we use static condensation to reduce the size of the global system to one involving just the pressure and trace degrees of freedom (all timings are for the entire solve, including the cost of static condensation); this reduces the cost of the smoother, at the expense of increasing the cost of prolongation operator construction, with a significant net savings in computational time and memory cost.

\begin{figure}[ht]
	\centering
		\scalebox{0.70}{
		\includegraphics{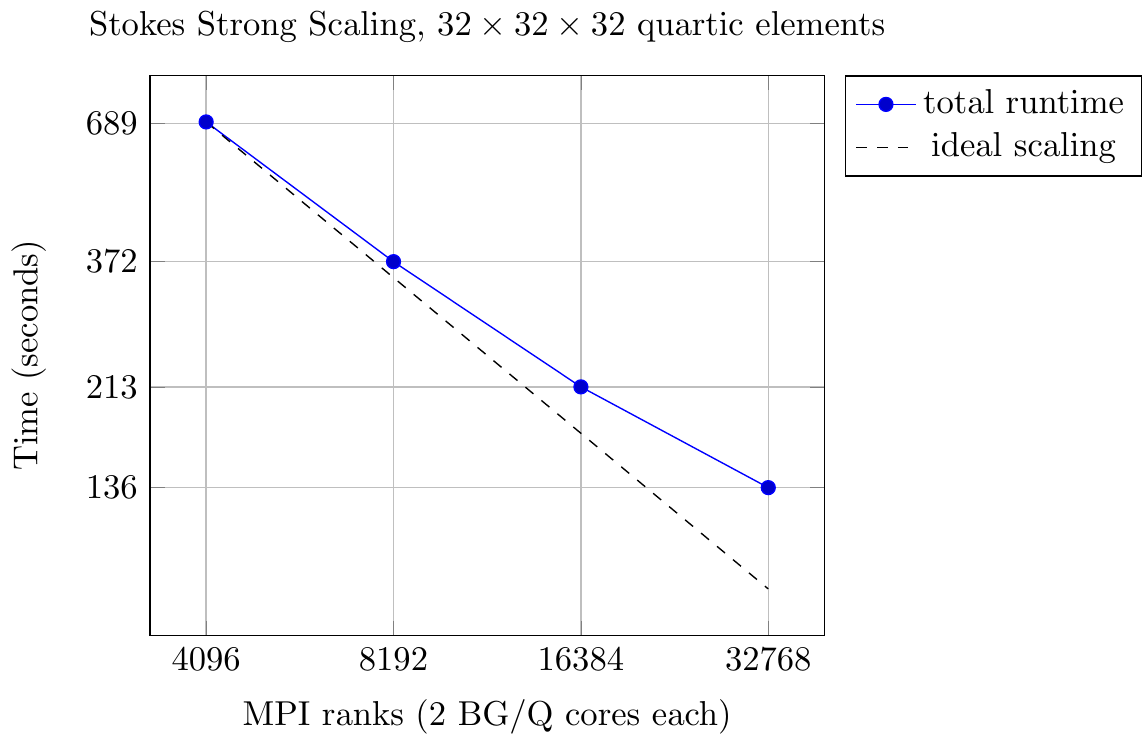}
		}
	\caption{Strong scaling for a Stokes solve with 76 million degrees of freedom, with 8 MPI ranks per node on Mira.}
\label{fig:runningTime}
\end{figure}

\begin{figure}[ht]
	\centering
	\scalebox{0.70}{
	\includegraphics{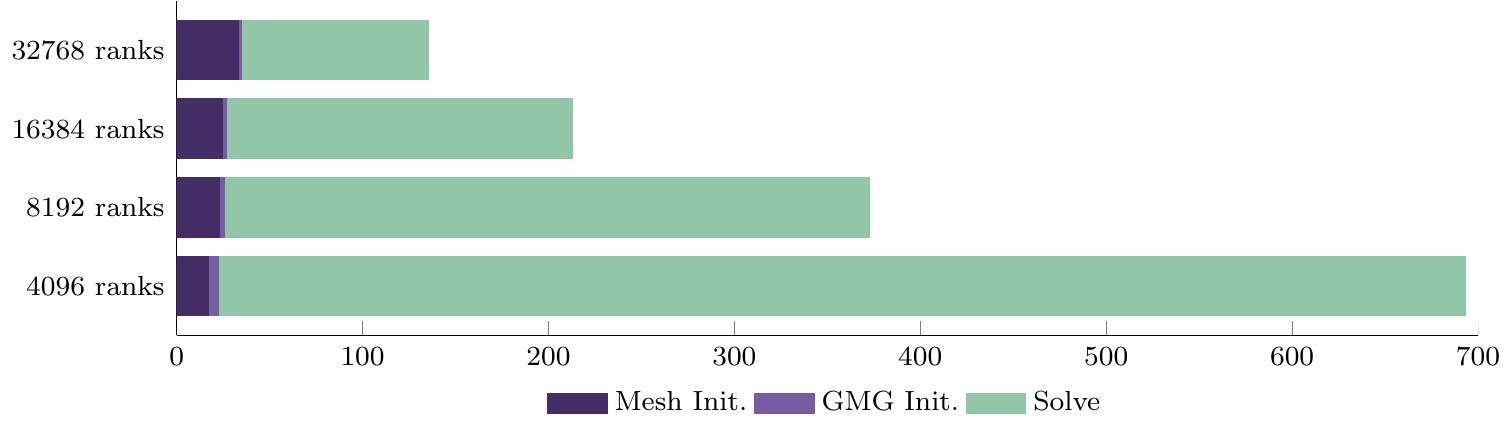}
	}
\caption{Timing detail for the solve whose strong scaling is shown in Figure \ref{fig:runningTime}.  ``GMG Init.'' refers to initialization of the geometric multigrid solver; the major cost within this is the labeling of the statically condensed degrees of freedom on the coarse meshes.  ``Solve'' includes determination of local stiffness contributions, as well as construction of the geometric multigrid prolongation and smoothing operators.  Further detail of the solve components can be seen in Figure \ref{fig:solveBreakdown}.}
\label{fig:timingBreakdown}
\end{figure}

The overall running times for our four runs are shown in Figure \ref{fig:runningTime}; the dashed line indicates ideal speedup starting from the smallest node count.  Even going to the limit of one element per rank, we observe a fivefold speedup (compared with the ideal eightfold speedup) in the overall runtime.  As can be seen in the detailed breakdown in Figure \ref{fig:timingBreakdown}, mesh initialization in particular shows no speedup.  We have recently introduced a distributed data structure for the mesh topology in Camellia---it may be that the lack of speedup is attributable to the communication costs associated with the distributed mesh; very likely with some further performance analysis and tuning the costs here could be reduced.  In the one-element-per-rank limit, mesh initialization takes about 25\% of the overall runtime.  If we neglect the cost of mesh initialization, then the scaling result is considerably better: the speedup is 6.5-fold.

\definecolor{mydeeppurple}{HTML}{442d65}
\definecolor{mypurple}{HTML}{775ba3}
\definecolor{mylightgreen}{HTML}{91c6a9}
\definecolor{myecru}{HTML}{F9E2B5}
\definecolor{myorange}{HTML}{FA8A5F}
\definecolor{mystudy}{HTML}{5F9EA2}

\begin{figure}[ht]
	\centering
	\scalebox{0.95}{
	\includegraphics{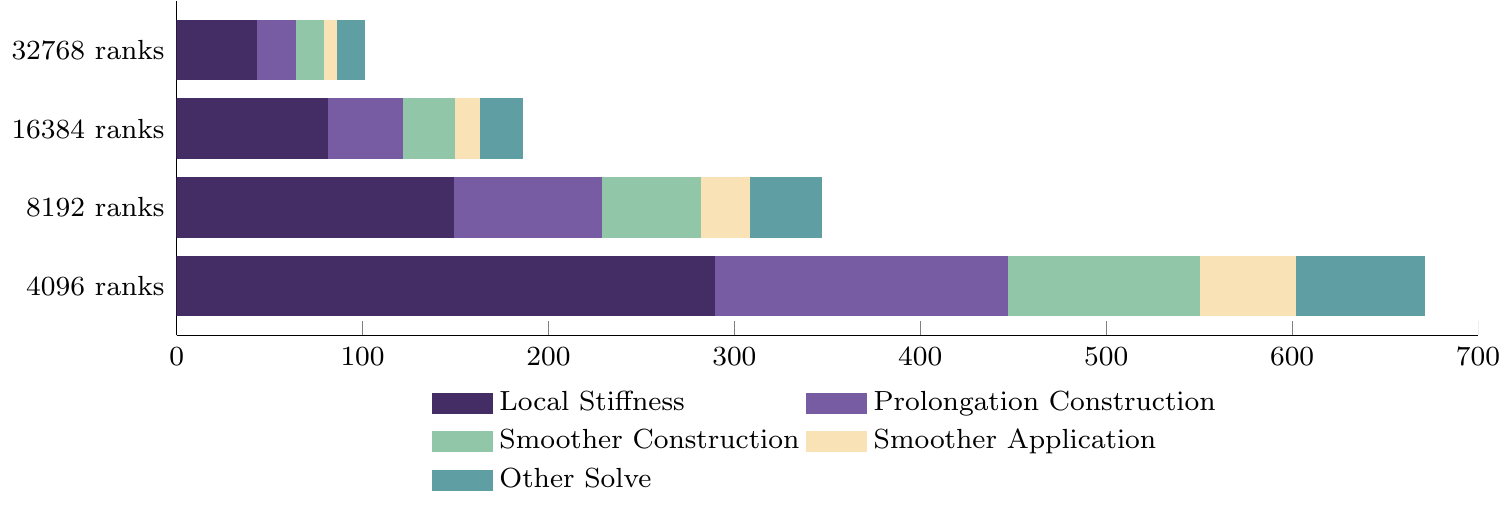}
	}
\caption{Breakdown of the components of the solve times shown shown in Figure \ref{fig:timingBreakdown}.  Every component of the solve demonstrates near-perfect scaling, up to the limit of one element per MPI rank.}
\label{fig:solveBreakdown}
\end{figure}

\subsection{Performance Enhancement Possibilities}
It is worth noting that, while we have optimized many parts of the code, there are still opportunities for substantially speeding things up and/or reducing the memory footprint further.  To name a few:
\begin{itemize}
\item For Stokes and Navier-Stokes, we could statically condense all but one of the pressure degrees of freedom on each mesh.
\item We could improve the Schwarz factorization by using a Cholesky factorization rather than the LU factorization; moreover, using Cholesky, we could reduce the memory footprint of the stored Schwarz factorizations by about half.
\item Similarly, we could improve our static condensation implementation by taking advantage of symmetry in the storage of element matrices and their factorization.
\end{itemize}

\pagebreak 
\section{Conclusion}
\label{sec:conclusion}
We have detailed a new approach to preconditioning DPG system matrices using geometric multigrid, and have demonstrated its efficiency for ultraweak Poisson, Stokes, and Navier-Stokes formulations through numerical experiments.  Moreover, ours is a \emph{black-box} approach, in that it can be applied to any DPG system (albeit with iterative performance characteristics that depend on the problem being solved).  We implemented our approach in Camellia, and have demonstrated the scalability of the implementation through experiments involving up to 32,768 MPI ranks.

\paragraph{License}The submitted manuscript has been created by UChicago Argonne, LLC, Operator of Argonne National Laboratory (``Argonne''). Argonne, a U.S. Department of Energy Office of Science laboratory, is operated under Contract No. DE- AC02-06CH11357. The U.S. Government retains for itself, and others acting on its behalf, a paid-up nonexclusive, irrevocable worldwide license in said article to reproduce, prepare derivative works, distribute copies to the public, and perform publicly and display publicly, by or on behalf of the Government.

\appendix

\clearpage

\section{Some Notes on the Implementation}
\label{sec:implementation}
While a full description of our implementation is beyond the scope of this paper (we hope to provide more detail in a forthcoming report on recent changes to Camellia), it is perhaps worth describing a couple of the core components at a high level.  Below, we describe a core mechanism by which the prolongation operator is constructed, then briefly describe our implementation of the overlapping additive Schwarz smoother.  To perform the conjugate gradient iteration, we make straightforward use of Trilinos's Belos package \cite{Amesos2Belos, Trilinos}.

\paragraph{Prolongation} For the construction of the prolongation operator, a crucial feature is the ability to represent a coarse basis in terms of a fine basis, for both $h$- and $p$-refinements of the coarse basis.  A standard approach would be to fix a family of basis functions, and to hard-code the relationships between the order $k$ and the order $\lfloor \frac{k}{2} \rfloor$ functions, as well as the relationships between the order $k$ basis on the coarse element and the bases on its $h$-refined children.  The advantage of this approach is that it will be computationally efficient.  The disadvantages of the approach are that it will require reimplementation for each family of basis functions and is potentially error-prone and difficult to debug.

Because we want to be flexible with regard to basis functions, we adopt an alternative approach, in which the relationship between coarse and fine bases is determined on the fly, through Camellia's \texttt{BasisReconciliation} class.  \texttt{BasisReconciliation} maintains a cache of previously computed relationships; because there are a limited number of these, the computational cost of computing the relationships is in practice negligible.  As arguments, the relevant \texttt{BasisReconciliation} methods take the bases on the coarse and fine elements and a \texttt{RefinementBranch} argument that specifies the geometric relationship between the coarse and fine elements.  (Though in all results here we take only a single refinement between mesh levels, our implementation supports an arbitrary number of refinements; it may be that performance improvements would be possible in some contexts by skipping some mesh levels.)

\paragraph{Smoothing} Trilinos's Ifpack package \cite{ifpack-guide} provides a host of features for incomplete matrix factorizations.  Included in Ifpack is an additive overlapping Schwarz implementation, \texttt{Ifpack\_AdditiveSchwarz}.  This implements a Schwarz factorization of an \texttt{Epetra} distributed sparse matrix at a specified level of overlap.  The blocks here are defined in terms of the distribution of the matrix; the overlap is also interpreted algebraically.

For our present purposes, we require finer control over the definition of the Schwarz blocks; we wish to define these geometrically in terms of degrees of freedom seen by a specific element or set of elements.  For this reason, we duplicated \texttt{Ifpack\_AdditiveSchwarz} to create an \texttt{AdditiveSchwarz} class within the Camellia namespace, and tailored it to provide these features---the \texttt{AdditiveSchwarz} constructor now takes as additional arguments a Camellia \texttt{Mesh} and \texttt{DofInterpreter}, which provide element information and define which degrees of freedom belong to each element.

Much of our implementation is the same as or similar to the original Ifpack implementation; among other things, this means that we have flexibility with regard to how the Schwarz blocks are factored.  In everything reported here, the factorization is done directly through the Amesos KLU solver; however, there is also support for e.g. incomplete Cholesky factorizations, and using these may provide some speedup relative to the results we have reported here.

\section{Schwarz Smoother Weights}
\label{sec:schwarzWeights}
In Section \ref{sec:sigmaWeightDiscussion} above, we describe two approaches to determining the scalar weight $\sigma$ that we apply to our smoother $B$.  The goal in selecting this weight is to keep the maximum eigenvalue of $\sigma B A$ at or below 1, to ensure the positive definiteness of the preconditioner; however, subject to that constraint, larger values of $\sigma$ will generally perform better.  The first approach uses $\sigma = \frac{1}{\Delta(G) + 2}$, where $\Delta(G)$ is the maximum degree of the adjacency graph $G$---with this choice, we can prove the eigenvalue bound holds.  However, the estimate may not be sharp, and in our present work we adopt a second approach that is generally more aggressive in its choice of $\sigma$.  In this approach, we count the face neighbors of the Schwarz overlap domain (including the elements that lie in the domain), and call the maximum such count $N$.  We then use $\sigma = \frac{1}{N+1}$.

For reference, the $\sigma$ values for uniform grids with overlap level $0$ (corresponding to our $p$-multigrid operators) are shown in Table \ref{table:SigmaValuesUsedOverlap0}.  The values for overlap level $1$ ($h$-multigrid) are shown in Table \ref{table:SigmaValuesUsedOverlap1}.  In every case, the $\frac{1}{N+1}$ value we use is at least as large as the alternative.

\begin{table}[h!b!p!]
\begin{center}
\begin{tabular}{ c  c  c  c  c}
\hline
Space Dimension 	& Mesh Width 		&$\frac{1}{\Delta(G) + 2}$		&$\frac{1}{N+1}$\\		
1				&2				&$1/4$					&$1/3$\\
2				&2				&$1/6$					&$1/4$\\
3				&2				&$1/10$					&$1/5$\\
\hline
1				&$>2$			&$1/4$					&$1/4$\\
2				&$>2$			&$1/6$					&$1/6$\\
3				&$>2$			&$1/10$					&$1/8$\\
\hline
\end{tabular}
\end{center} 
\caption[]{Two possible choices for $\sigma$ on uniform grids, when overlap level is $0$ (used for $p$-multigrid operators).  We use the more aggressive $\frac{1}{N+1}$ choice in our experiments.}
\label{table:SigmaValuesUsedOverlap0}
\end{table}

\begin{table}[h!b!p!]
\begin{center}
\begin{tabular}{ c  c  c  c  c}
\hline
Space Dimension 	& Mesh Width 		&$\frac{1}{\Delta(G) + 2}$		&$\frac{1}{N+1}$\\		
1				&2				&$1/4$					&$1/3$\\
2				&2				&$1/6$					&$1/5$\\
3				&2				&$1/10$					&$1/8$\\
\hline
1				&4				&$1/6$					&$1/5$\\
2				&4				&$1/14$					&$1/12$\\
3				&4				&$1/34$					&$1/23$\\
\hline
1				&$>4$			&$1/6$					&$1/6$\\
2				&$>4$			&$1/14$					&$1/14$\\
3				&$>4$			&$1/34$					&$1/26$\\
\hline
\end{tabular}
\end{center} 
\caption[]{Two possible choices for $\sigma$ on uniform grids, when overlap level is $1$ (used for $h$-multigrid operators).  We use the more aggressive $\frac{1}{N+1}$ choice in our experiments.}
\label{table:SigmaValuesUsedOverlap1}
\end{table}

\clearpage

\bibliographystyle{plain}

\begin{thebibliography}{10}

\bibitem{achdou1999iterative}
Yves Achdou, Yvon Maday, and Olof Widlund.
\newblock Iterative substructuring preconditioners for mortar element methods
  in two dimensions.
\newblock {\em SIAM Journal on Numerical Analysis}, 36(2):551--580, 1999.

\bibitem{Arioli2003}
M.~Arioli.
\newblock A stopping criterion for the conjugate gradient algorithm in a finite
  element method framework.
\newblock {\em Numerische Mathematik}, 97(1):1--24, 2003.

\bibitem{BarkerDPG}
Andrew~T Barker, Susanne~C Brenner, Eun-Hee Park, and Li-Yeng Sung.
\newblock A one-level additive {S}chwarz preconditioner for a discontinuous
  {P}etrov-{G}alerkin method.
\newblock In {\em Domain Decomposition Methods in Science and Engineering XXI},
  pages 417--425. Springer, 2014.

\bibitem{Amesos2Belos}
Eric Bavier, Mark Hoemmen, Sivasankaran Rajamanickam, and Heidi Thornquist.
\newblock Direct and iterative solvers for large sparse linear systems.
\newblock {\em Scientific Programming}, 20(3), 2012.

\bibitem{bey1997downwind}
J{\"u}rgen Bey and Gabriel Wittum.
\newblock Downwind numbering: Robust multigrid for convection-diffusion
  problems.
\newblock {\em Applied Numerical Mathematics}, 23(1):177--192, 1997.

\bibitem{bochev2009least}
Pavel~B Bochev and Max~D Gunzburger.
\newblock {\em Least-squares finite element methods}, volume 166.
\newblock Springer Science \& Business Media, 2009.

\bibitem{Intrepid}
Pavel~B. Bochev, Robert~C. Kirby, Kara~J. Peterson, and Denis Ridzal.
\newblock Intrepid project.
\newblock {\em http://trilinos.sandia.gov/packages/intrepid/}.

\bibitem{BramwellDemkowiczGopalakrishnanQiu11}
J.~Bramwell, L.~Demkowicz, J.~Gopalakrishnan, and W.~Qiu.
\newblock A locking-free $hp$ {DPG} method for linear elasticity with symmetric
  stresses.
\newblock {\em Num. Math.}, 122(4):671--707, December 2012.

\bibitem{BrennerScott}
Susanne~C. Brenner and L.~Ridgeway Scott.
\newblock {\em The Mathematical Theory of Finite Element Methods}.
\newblock Springer, 3rd edition, 2008.

\bibitem{Brooks1941}
R.~L. Brooks.
\newblock On colouring the nodes of a network.
\newblock {\em Mathematical Proceedings of the Cambridge Philosophical
  Society}, 37:194--197, 4 1941.

\bibitem{cai1997first}
Zhiqiang Cai, Thomas~A Manteuffel, and Stephen~F McCormick.
\newblock First-order system least squares for second-order partial
  differential equations: Part {II}.
\newblock {\em SIAM Journal on Numerical Analysis}, 34(2):425--454, 1997.

\bibitem{DPGCompressible}
Jesse Chan, Leszek Demkowicz, and Robert Moser.
\newblock A {DPG} method for steady viscous compressible flow.
\newblock {\em Computers \& Fluids}, 98:69 -- 90, 2014.
\newblock 12th {USNCCM} mini-symposium of High-Order Methods for Computational
  Fluid Dynamics - A special issue dedicated to the 80th birthday of Professor
  Antony Jameson.

\bibitem{CockburnKanschatSchotzauSchwab03}
B.~Cockburn, G.~Kanschat, D.~Schotzau, and Ch. Schwab.
\newblock Local {D}iscontinuous {G}alerkin methods for the {S}tokes system.
\newblock {\em SIAM J. on Num. Anal.}, 40:319--343, 2003.

\bibitem{dahmen2012adaptive}
Wolfgang Dahmen, Chunyan Huang, Christoph Schwab, and Gerrit Welper.
\newblock Adaptive {Petrov--Galerkin} methods for first order transport
  equations.
\newblock {\em SIAM journal on numerical analysis}, 50(5):2420--2445, 2012.

\bibitem{DPG1}
L.~Demkowicz and J.~Gopalakrishnan.
\newblock A class of discontinuous {Petrov-Galerkin} methods. {Part I}: {T}he
  transport equation.
\newblock {\em Comput. Methods Appl. Mech. Engrg.}, 199:1558--1572, 2010.
\newblock See also ICES Report 2009-12.

\bibitem{DPG6}
L.~Demkowicz and J.~Gopalakrishnan.
\newblock Analysis of the {DPG} method for the {P}oisson problem.
\newblock {\em SIAM J. Num. Anal.}, 49(5):1788--1809, 2011.

\bibitem{DPG2}
L.~Demkowicz and J.~Gopalakrishnan.
\newblock A class of discontinuous {Petrov-Galerkin} methods. {Part II}:
  {O}ptimal test functions.
\newblock {\em Numer. Meth. Part. D. E.}, 27(1):70--105, January 2011.

\bibitem{demkowicz2015ices}
Leszek Demkowicz and Jay Gopalakrishnan.
\newblock Discontinuous {Petrov-Galerkin (DPG)} method.
\newblock {\em ICES Report}, (15-20), 2015.

\bibitem{Fischer2005}
Paul~F. Fischer and James~W. Lottes.
\newblock {\em Hybrid Schwarz-Multigrid Methods for the Spectral Element
  Method: Extensions to Navier-Stokes}, pages 35--49.
\newblock Springer Berlin Heidelberg, Berlin, Heidelberg, 2005.

\bibitem{GopalakrishnanQiu11}
Jay Gopalakrishnan and Weifeng Qiu.
\newblock {An analysis of the practical DPG method}.
\newblock {\em Mathematics of Computation}, 2012.

\bibitem{gopalakrishnan2015degree}
Jay Gopalakrishnan and Joachim Sch{\"o}berl.
\newblock Degree and wavenumber [in] dependence of {S}chwarz preconditioner for
  the {DPG} method.
\newblock In {\em Spectral and High Order Methods for Partial Differential
  Equations ICOSAHOM 2014}, pages 257--265. Springer, 2015.

\bibitem{Trilinos}
Michael~A. Heroux, Roscoe~A. Bartlett, Vicki~E. Howle, Robert~J. Hoekstra,
  Jonathan~J. Hu, Tamara~G. Kolda, Richard~B. Lehoucq, Kevin~R. Long, Roger~P.
  Pawlowski, Eric~T. Phipps, Andrew~G. Salinger, Heidi~K. Thornquist, Ray~S.
  Tuminaro, James~M. Willenbring, Alan Williams, and Kendall~S. Stanley.
\newblock {An overview of the Trilinos project}.
\newblock {\em ACM Trans. Math. Softw.}, 31(3):397--423, 2005.

\bibitem{Kovasznay}
L.~I.~G. Kovasznay.
\newblock Laminar flow behind a two-dimensional grid.
\newblock {\em Mathematical Proceedings of the Cambridge Philosophical
  Society}, 44(01):58--62, 1948.

\bibitem{loghin1997preconditioning}
Daniel Loghin and Andrew~J Wathen.
\newblock Preconditioning the advection-diffusion equation: the green's
  function approach.
\newblock 1997.

\bibitem{RobertsDissertation}
Nathan~V. Roberts.
\newblock {\em A Discontinuous {Petrov-Galerkin} Methodology for Incompressible
  Flow Problems}.
\newblock PhD thesis, University of Texas at Austin, 2013.

\bibitem{Camellia}
Nathan~V. Roberts.
\newblock Camellia: A software framework for discontinuous {P}etrov-{G}alerkin
  methods.
\newblock {\em Computers \& Mathematics with Applications}, 2014.

\bibitem{DPGStokes}
Nathan~V. Roberts, Tan Bui-Thanh, and Leszek~F. Demkowicz.
\newblock The {DPG} method for the {S}tokes problem.
\newblock {\em Computers and Mathematics with Applications}, 2014.

\bibitem{DPGNavierStokes}
Nathan~V. Roberts, Leszek Demkowicz, and Robert Moser.
\newblock A discontinuous {P}etrov--{G}alerkin methodology for adaptive
  solutions to the incompressible {N}avier--{S}tokes equations.
\newblock {\em Journal of Computational Physics}, 301:456 -- 483, 2015.

\bibitem{ifpack-guide}
M.~Sala and M.~Heroux.
\newblock Robust algebraic preconditioners with {IFPACK} 3.0.
\newblock Technical Report SAND-0662, Sandia National Laboratories, 2005.

\bibitem{SmithBjorstadGropp}
Barry Smith, Petter Bj{\o}stad, and William Gropp.
\newblock {\em Domain Decomposition}.
\newblock Cambridge University Press, 1996.

\bibitem{toselli2005domain}
Andrea Toselli and Olof~B Widlund.
\newblock {\em Domain decomposition methods: algorithms and theory}, volume~34.
\newblock Springer, 2005.

\bibitem{WienersWohlmuth}
C.~Wieners and B.~Wohlmuth.
\newblock {Robust operator estimates and the application to sub structuring
  methods for first-order systems}.
\newblock {\em ESAIM: Mathematical Modelling and Numerical Analysis},
  48(5):1473--1494, 2014.

\bibitem{xu2016domain}
Xuejun Xu and Xiang Li.
\newblock Domain decomposition preconditioners for the discontinuous
  {Petrov-Galerkin} method.
\newblock {\em {ESAIM: Mathematical Modelling and Numerical Analysis}}, 2016.

\end{thebibliography}

\end{document}